\documentclass[11pt]{amsart}

\usepackage[T1]{fontenc}
\usepackage[utf8]{inputenc}
\usepackage[english]{babel}

\usepackage{amsmath,amssymb,amsfonts,amsthm,amsxtra,latexsym,mathtools,esint,tikz}

\usepackage{ifthen}
\usepackage{enumitem}
\usepackage{graphicx}
\usepackage{tikz}
\usepackage{xcolor}
\usepackage{soul}
\usepackage{cite}
\usepackage{microtype}
\usepackage{appendix}
\usepackage[left=2.5cm,right=2.5cm,top=2.7cm,bottom=2.7cm]{geometry}

\usepackage[colorlinks=true,urlcolor=blue,citecolor=red,linkcolor=blue]{hyperref}

\numberwithin{equation}{section}
\newtheorem{theorem}{Theorem}[section]
\newtheorem{lemma}{Lemma}[section]

\newtheorem{corollary}{Corollary}[section]
\newtheorem{propo}{Proposition}[section]

\DeclareMathOperator*{\osc}{osc}

\newcommand{\N}{{\mathbb N}}

\newcommand{\R}{{\mathbb R}}
\newcommand{\A}{\mathbb{A}}
\newcommand{\B}{\mathcal{B}}
\newcommand{\C}{\mathcal{C}}
\newcommand{\E}{\mathcal{E}}
\def\H{\mathcal{H}}
\newcommand{\I}{\mathcal{I}}
\newcommand{\K}{\mathcal{K}}
\def\L{\mathcal{L}}
\def\P{\mathcal{P}}
\newcommand{\Q}{\mathcal{Q}}
\def\S{\mathcal{S}}
\newcommand{\T}{\mathcal{T}}

\def\b{\mathbb{B}}
\def\k{\mathbb{K}}
\def\q{\mathbb{Q}}

\newcommand{\eps}{\varepsilon}
\newcommand{\e}{\epsilon}
\newcommand{\tr}{(u-k)_{-}}
\newcommand{\dive}{\mathrm{div}}
\def\d{\mathrm{d}}
\def\dim{\mathrm{dim}}
\newcommand{\dist}{\mathrm{dist}}

\newcommand{\beq}{\begin{equation}}
\newcommand{\eeq}{\end{equation}}

\def\XXint#1#2#3{\setbox0=\hbox{$#1{#2#3}{\int}$}
   \vcenter{\hbox{$#2#3$}}\kern-.5\wd0}

\newcommand{\llbracket}{\mathopen{\lbrack\!\lbrack}}
\newcommand{\rrbracket}{\mathclose{\rbrack\!\rbrack}}

\title[Anisotropic Trudinger's Equation]{%
Qualitative properties of solutions to parabolic anisotropic equations,
Part II.\vskip0.2cm
\noindent
- The anisotropic Trudinger's equation-
}
\author[Ciani, Henriques, Savchenko, Skrypnik]{S. Ciani {\&} E. Henriques {\&} M. O. Savchenko {\&} I. I. Skrypnik}
\address{Department of Mathematics of the University of Bologna, Piazza Porta San Donato, 5, 40126 Bologna, Italy}
\email{simone.ciani3@unibo.it}
\address{Centro de Matem\'atica, Universidade do Minho - Polo CMAT-UTAD Departamento de Matem\'atica Universidade de Tr\'as-os-Montes e Alto Douro, Vila Real, Portugal}
\email{eurica@utad.pt}
\address{Institute of Applied Mathematics and Mechanics, National Academy of Sciences of Ukraine, Gen. Batiouk Str. 19, 84116 Sloviansk, Ukraine}
\email{ihor.skrypnik@gmail.com {\&} shan\_maria@ukr.net}

\begin{document}
\begin{abstract} 
We study the local regularity properties of weak solutions to a special class of anisotropic doubly nonlinear parabolic operators, whose prototype is the anisotropic Trudinger's equation
$$
u_t- \sum\limits_{i=1}^N D_i\Big(u^{2-p_i}|D_i u|^{p_i-2} D_i u\Big)=0,\qquad u\geqslant 0.
$$
We prove a parabolic Harnack inequality, valid without any restrictions on the exponents $p_i$s. When the range of diffusion exponents is restricted, solutions are H\"{o}lder continuous. \vskip0.2cm \noindent 

\end{abstract}

\maketitle
	\begin{center}
		\begin{minipage}{9cm}
			\small
   \tableofcontents
		\end{minipage}
	\end{center}



\def\tr{(u-k)_{-}}
\def\R{\mathbb{R}}
\def\A{\mathbb{A}}
\def\Q{\mathcal{Q}}
\def\b{\mathbb{B}}
\def\K{\mathcal{K}}
\def\k{\mathbb{K}}
\def\E{\mathcal{E}}
\def\q{\mathbb{Q}}

\def\N{\mathbb{N}}
\def\dive{\mathrm{div}}
\def\d{\mathrm{d}}
 \def\B{\mathcal{B}}
\def\dim{\mathrm{dim}}
\def\dist{\mathrm{dist}}
\def\H{\mathcal{H}}
\def\I{\mathcal{I}}
\def\L{\mathcal{L}}
\def\C{\mathcal{C}}
\def\T{\mathcal{T}}
\def\S{\mathcal{S}}
\def\P{\mathcal{P}}
\def\e{\epsilon}

\addtocontents{toc}{\protect\setcounter{tocdepth}{1}}

\section{Introduction and Main results}\label{Sec.1}
\subsection{Introduction}
In the framework of our program on the qualitative properties of local weak solutions to parabolic anisotropic operators, here our motivation stems from the classic research in \cite{Tru} and  the relatively recent studies in \cite{KinKuu, KuuSilUrb, KuuLalSilUrb, GiaVes, Lia}, where the regularity results such as the local boundedness, Harnack's inequality, H\"{o}lder continuity of solutions to doubly nonlinear parabolic equations
$$u_t -div\Big(u^{2-p}|D u|^{p-2} D u\Big)=0,\quad u\geqslant 0,$$
were proved.  Such type of equations in the literature are known as Trudinger's equations, following his pioneering work in \cite{Tru}.  In this paper we are concerned with  non-negative solutions to anisotropic Trudinger's equations with measurable and bounded coefficients, whose prototype is
\begin{equation} \label{prototype} u_t-\sum\limits_{i=1}^N D_i\Big(u^{2-p_i}|D_i u|^{p_i-2} D_i u\Big)=0,\quad u\geqslant 0.\end{equation} We observe that when $p_i \equiv p$, equation \eqref{prototype} enjoys a similar structure to Trudinger's equation, with the only difference that the modulus of ellipticity depends on each directional derivatives. In this case, the change of variables $v=u^{\frac{1}{p-1}}$, turns  equation \eqref{prototype} into
\begin{equation}
    \label{proto-2} 
    \Big( v^{p-1}\Big)_t + \sum_{i=1}^N \gamma_i D_i \, \Big(|D_i v|^{p-2} D_i v \Big) =0\,, \qquad \gamma_i >0\,,
\end{equation} whose parabolic and elliptic energy terms are homogeneous with respect to the multiplication by a constant. The original idea of \cite{Tru} was that this homogeneity makes the low-order regularity theory simpler, but on the other hand the treatment of the time derivative is encumbered, at the point that adding or subtracting a constant to a solution does not produce a solution to the same kind of equation. For this reason, the study of the Harnack inequality and of the H\"older continuity seems to be dinstinct, see for instance \cite{BDL21} and  \cite{GiaVes}. Here below we address both studies, for operators whose prototipe is \eqref{prototype}.\vskip0.1cm \noindent More precisely, let $\Omega\subset \R^N$ be open and bounded, $T>0$ and let us consider  numbers 
\[ 1<p_1 \leq p_2 \leq \dots \leq p_N < \infty \, .\]
\noindent We study the local regularity properties of non-negative solutions to the equation
\begin{equation}\label{eq1.1}
u_{t}-\sum\limits_{i=1}^N D_i A_i (x, t, u, D u)=0, \quad \text{locally weakly in} \quad  \Omega_{T}:=\Omega \times (0, T).
\end{equation}
where $A_i(x,t,s,\xi):\Omega_{T}\times \mathbb{R}_+\times \mathbb{R}^{N} \rightarrow \mathbb{R}^{N}$ are Carath\'eodory functions, satisfying the following structure conditions
\begin{equation}\label{eq1.2}
\begin{cases}
 A_i(x, t, s, \xi)\, \xi_i  \geqslant  K_{1}\, s^{2-p_i}\,|\xi_i|^{p_i},\qquad \forall (s,\xi) \in \R_+ \times \R^N,\\ $\quad$ \\
|A_i(x, t, s,  \xi)| \leqslant  K_{2}\, s^{2-p_i} |\xi_i|^{p_i-1},\qquad  \forall \, i=1, ...,N,
\end{cases}
\end{equation}
for positive constants $K_1,K_2$. The formulation of the structure conditions \eqref{eq1.2} can be rephrased under the assumption of the chain rule, see Section \ref{redefinition} for the details.

\subsection{Main Results} \label{sec-main results}
Our first result is the Harnack inequality. The special feature of the following theorem being its validity regardless of any restrictions on the gap between $p_N$ and $p_{1}$. 

\noindent The parameters $N$, $p_1$, ..., $p_N$, $K_1$, $K_2$ will be referred to as {\it the data}. 

\begin{theorem}\label{th1.1}
Let $u$ be a non-negative, local weak solution to \eqref{eq1.1}-\eqref{eq1.2} in $\Omega_T$. There exist positive constants $C$, $\bar{C}$  depending only on the data such that, for all cylinders \[
\K_{8\rho}(x_0) \times (t_0- \bar{C}(8\rho)^{p}, t_0+ \bar{C}(8\rho)^{p}) \subset \Omega_T \ , \] there holds
\begin{equation}\label{eq1.5}
\frac{1}{C}\sup\limits_{\K_\rho(x_0)}u(\cdot, t_0-\bar{C}\,\rho^{p})\leqslant u(x_0, t_0)\leqslant C\, \inf\limits_{\K_\rho(x_0)} u(\cdot, t_0+\bar{C}\,\rho^{p}),
\end{equation} 
where
\[
\K_{\rho}(x_0):=\prod_{i=1}^N \bigg\{|x_i-{x_0}_i|<\rho^{\frac{p}{p_i}}\bigg\}\,\quad \text{and} \qquad \frac{1}{p}:=\frac{1}{N}\sum\limits_{i=1}^N\frac{1}{p_i} \,.\]
\end{theorem}
\noindent Our second result is the local H\"{o}lder continuity of solutions. The proof distinguishes between the cases when $u$ is large or close to zero. In the former case, the equation can be regarded as the anisotropic $p$-Laplacian (cf. \cite{CHS24}, \cite{CiaMosVes}), whose regularity theory is known only under the additional assumption that the gap between $p_N$ and $p_{1}$ is sufficiently small (see Lemma \ref{lem2.8}). For this reason a constraint on the spareness of the $p_i$s appears. See also \cite{DiBGiaVes3} for this technique in the case of porous medium equations.

\begin{theorem}\label{th1.2}
Let $u$ be a non-negative, local weak solution to \eqref{eq1.1}-\eqref{eq1.2}. There exists a number $\epsilon_*\in (0, 1)$ depending only on the data such that, if 
\begin{equation}\label{eq1.6}
p_N-p_{1}\leqslant \epsilon_*,
\end{equation} then $u$ has a locally H\"{o}lder continuous representative.
\end{theorem}

\subsection{Novelty and Significance} The issues touched in here have been considered in the standard case $p_i \equiv p$ starting with  Caffarelli and Friedman \cite{CafFri}, Caffarelli and Evans
\cite{CafEva},  DiBenedetto and Friedman \cite{DiBFri} and DiBenedetto \cite{DiB} and continued in \cite{PorVes, Ves1, Ves2, Iva1, Iva2, BogDuzGiaLiaSch, DiBGiaVes1, DiBGiaVes2, DiBGiaVes3}.  Since Trudinger's equation seems to be one of the few nonlinear parabolic equation which admits  a
Harnack inequality just as for the heat equation\footnote{Meaning with this that it does not require to consider an intrinsic geometry (see for this \cite{DiB}).}, it is the natural candidate for generalizing these results
to nonlinear equations with nonstandard growth, particularly to anisotropic equations. The study of regularity properties of solutions to parabolic equations with nonstandard growth goes back to the 1970's. Kolodii \cite{Kol1, Kol2} gave his seminal contribution on the local boundedness after which several other came along,  see for example \cite{MinXit, Sin, Str}. At the same time, the basic qualitative properties such as  continuity and the Harnack inequality of solutions to anisotropic elliptic and parabolic equations are still poorly understood.  Regarding the continuity of solutions, surprisingly, the most extensively studied case
is the anisotropic porous medium equation \cite{BurSkr, HU06, Hen1, Hen2, EH21}. To the best of our knowledge, partial results for anisotropic elliptic equations can be found in \cite{LiaSkrVes, CiaSkrVes, LisSkr, CiaHenSkr} or in \cite{DiBGiaVes4} under the additional condition that the gap between $p_N$ and $p_{1}$ is small enough. For the elliptic prototype, the authors in \cite{Brasco} show that locally bounded weak solutions are Lipschitz continuous, and hence continuous, without any restrictions on the spread of the $p_i$s. As for the anisotropic parabolic equations, some Harnack-type inequalities were obtained \cite{CiaGua, CiaMosVes, CHS24, CiaHen}. Questions on boundedness in a wider anisotropic regime were addressed in \cite{CiaHenSavSkr1}, while in \cite{CiaHenSavSkr2} the authors proved the expansion of positivity. In this paper we continue the study of the regularity properties for anisotropic parabolic equations started in \cite{CiaHenSavSkr1, CiaHenSavSkr2}, namely we prove Harnack's inequality and  H\"{o}lder continuity of weak solutions to a class of anisotropic parabolic equations having as prototype the  Trudinger's equation \eqref{eq1.1}. In this context, it is a striking fact that the boundedness (proved in \cite{CiaHenSavSkr1}) and Harnack's inequality for solutions of such equations can be proved without any restrictions on the numbers $p_i$ and have the same form as for the heat equation (cf. Theorem \ref{th1.1}).

\subsubsection{The Method} \noindent In \cite{DeG} Ennio De Giorgi developed an original geometric method for the boundedness and regularity
of solutions to elliptic equations with discontinuous, but measurable and bounded coefficients. The fundamental ideas of this
technique have been successfully applied to get regularity for local minima of calculus of
variations and/or solutions to the corresponding elliptic and parabolic equations with standard
growth (for an exhaustive overview on the subject, see, for example \cite{DiB, LadUra, LadSolUra, DiBGiaVes3}). The core ideas of our proofs are based on this method and also on Krylov and Safonov \cite{KrySaf} procedure, when proving the Harnack inequality, and on the intrinsic scaling method (developed by DiBenedetto \cite{DiB}) when proving the H\"{o}lder continuity of solutions. 
Although the strategy for the proof of our main results goes as in the standard $p$-growth, we had to overcome the difficulties inherent to the presence of the anisotropy. The proof of the H\"older continuity follows the lines of \cite{KuuSilUrb} and \cite{KuuLalSilUrb}, with the complication that the case where $u$ is away from zero does not allow us to transform \eqref{prototype} into the $p$-Laplacian operator, for which there would have been already a complete picture. Instead, it turns \eqref{prototype} into the anisotropic $p_i$-Laplacian, for which the regularity theory is still under development. 

\subsubsection{Structure of the paper} The paper is organized as follows. In Section $2$ we collect some preliminary and technical
properties, being Section $3$ devoted to the proof of a De Giorgi type lemma. Expansion of positivity is proved in Section $4$. Our first main result, the Harnack inequality, is proved in Section $5$. Finally, in Section $6$, we prove H\"{o}lder continuity
of weak solutions to \eqref{eq1.1}. Motivation for the functional framework, the definition of solution with the exponential mollification and other tools of the trade are collected in the Appendix \ref{appendix}.

\section{Notation}
\noindent \subsection{Indexes} From now on we assume that the numbers $\{p_1,p_2, \dots, p_N\}$ are fixed. In order to ease notation we further assume, without loss of generality, that 
\[1<p_1 \leq p_2 \leq \dots \leq p_N< \infty\] and we set
\[\frac{1}{p}:=\frac{1}{N}\sum\limits_{i=1}^N\frac{1}{p_i}= \frac{1}{N}\sum\limits_{i}\frac{1}{p_i}\,.\]
When index $i$ is concerned, we will always reduce the set of indexes as in the previous formula, also for the product $\prod_{i=1}^N= \prod_i$. Given $K_1,K_2$ by the structure conditions \eqref{eq1.2}, we will refer to the set of parameters $\{N, p_1, \dots, p_N, K_1,K_1\}$ as the set of (structural) {{\it data}, and by saying that a constant $\gamma$ depends only on the data we mean that it can be quantitatively determined {\it a priori} only in terms of the above quantities. As usual, the generic constant $\gamma$ may change from line to line. 
\subsection{Geometry} For fixed $r, \eta >0$ and a point $(\bar{x},\bar{t}) \in \R^{N+1}$, we define the standard cube
\[K_r(\bar{x})= \big\{x\in \R^N:\quad |x_i-\bar{x}_i|<r, \quad i=1, \dots, N \big\} \ , \]
the anisotropic cubes
\[\K_{r}(\bar{x}):=\big\{x\in \R^N:\quad |x_i-\bar{x}_i|<r^{\frac{p}{p_i}},\quad i=1, ...,N\big\},\]
and the cylinders
\[Q_{r, \eta}(\bar{x}, \bar{t}):=\K_{r}(\bar{x})\times (\bar{t}-\eta, \bar{t}).\]
Finally, for fixed $k,r,\eta>0$, we define the intrinsic anisotropic cubes 
\[\K^k_r(\bar{x}):=\big\{x\in \R^N:\quad  |x_i-\bar{x}_i|<r^{\frac{p}{p_i}} k^{-\frac{p_N-p_i}{p_i}},\quad i=1, ..., N\big\} \ ,\]
and cylinders
\[Q^k_{r, \eta}(\bar{x}, \bar{t}):=\K^k_r(\bar{x})\times (\bar{t}-\eta, \bar{t}).\]

\section{Functional Setting and Auxiliary Results}\label{Sec.2}

\subsection{Functional Anisotropic Setting and Parabolic Embeddings}
\noindent Equation \eqref{prototype} has to be understood in an appropriate variational meaning. For this aim, we define the anisotropic spaces of locally integrable functions as
$$W^{1,\textbf{p}}(\Omega)=\big\{u\in W^{1,1}(\Omega): D_i u\in L^{p_i}(\Omega),\quad i=1, ...,N\big\},$$
$$L^{\textbf{p}}(0, T; W^{1,\textbf{p}}(\Omega))=\big\{u\in L^1(0, T; W^{1,1}(\Omega)): D_i u\in L^{p_i}(0, T; L^{p_i}(\Omega)),\quad i=1, ...,N\big\},$$
and the respective spaces of functions with zero boundary data
$$W^{1,\textbf{p}}_0(\Omega)=W^{1,1}_0(\Omega)\cap W^{1,\textbf{p}}(\Omega),$$
$$L^{\textbf{p}}(0, T; W^{1,\textbf{p}}_0(\Omega))=L^{1}(0, T; W^{1,1}_0(\Omega))\cap L^{\textbf{p}}(0, T; W^{1,\textbf{p}}(\Omega)),$$
for ${\bf{p}}=(p_1, \cdots, p_N)$.
The following lemma was originally studied in \cite{Tro}, see \cite{Mosconi} for a modern proof.
\begin{lemma}\label{lem2.1}
Let $\Omega\subset \R^N$ be open and bounded, $u\in W^{1,1}_0(\Omega)$ and let us assume
\[\sum_i \int\limits_{\Omega} |D_i |u|^{\alpha_i}|^{p_i} dx <\infty,\quad \text{for} \quad \alpha_i>0\quad \text{and} \quad p_i>0 \, \, \text{such that} \quad p<N\,.\]
Then there exists $c>0$, depending on $\{ N, p_1, \dots, p_N, \alpha_1, \dots, \alpha_N\}$, such that
\begin{equation}\label{eq2.1}
\Big(\int\limits_{\Omega}|u|^{p^*_{\alpha}} dx\Big)^{\frac{N-p}{N}}\leqslant c \prod_i \Big(\int\limits_{\Omega}|D_i |u|^{\alpha_i}|^{p_i} dx \Big)^{\frac{p}{N p_i}},\quad p^*_{\alpha}=\frac{N \alpha p}{N-p},\quad \alpha=\frac{1}{N}\sum_i \alpha_i.
\end{equation}
\end{lemma}
\noindent By using H\"{o}lder's inequality and \eqref{eq2.1}, the following lemma is valid too.
\begin{lemma}\label{lem2.2}
Let $\Omega\subset \R^N$ be open and
bounded, $u\in L^1(0, T; W^{1,1}_0(\Omega))$, $p<N$ and let also
\[\sup\limits_{0<t<T}\int\limits_{\Omega} |u|^\beta \,dx<\infty,\quad \beta>0\, ,\]
and
\[ \sum_i \int\limits^T_0\int\limits_{\Omega} |D_i |u|^{\alpha_i}|^{p_i} dx\,dt <\infty,\quad \alpha_i>0,\quad p_i>1,\quad i=1, ...,N.\]
Then there exists $c>0$, depending on $N, p_1, ..., p_N, \alpha_1, ..., \alpha_N$, such that
\begin{multline}\label{eq2.2}
\int\limits^T_0\int\limits_{\Omega}|u|^{p(\alpha+\frac{\beta}{N})} dx\,dt\leqslant c \Big(\sup\limits_{0<t<T}\int\limits_\Omega |u|^\beta\,dx\Big)^{\frac{p}{N}}\prod_i \int\limits^T_0\Big(\int\limits_{\Omega}|D_i |u|^{\alpha_i}|^{p_i} dx \Big)^{\frac{p}{N p_i}}\,dt\\\leqslant c \Big(\sup\limits_{0<t<T}\int\limits_\Omega |u|^\beta\,dx\Big)^{\frac{p}{N}}\bigg[\sum_i \int\limits^T_0\int\limits_{\Omega}|D_i |u|^{\alpha_i}|^{p_i} dx\,dt\bigg].
\end{multline}
\end{lemma}
\noindent A closer inspection of our proofs (Lemma \ref{lem3.1} here below, and Lemma \ref{lem2.7} of \cite{CiaHenSavSkr1}) reveals, as in Remark 3.9 of \cite{Vestberg}, that the condition $p<N$ can be removed.

\subsection{Definition of weak solution} \label{def-sol} Here below we define the notion of local weak solution to
equation \eqref{eq1.1}. For the motivation lying behind this definition and an alternative definition that takes into account an approximated time-derivative, we refer to section \ref{mollificazione}. A measurable function $u : \Omega_T \rightarrow [0, \infty)$ such that
\begin{equation} \label{functional}
\begin{cases} 
u\in C_{loc}(0, T; L^{\frac{p_N}{p_N-1}}_{loc}(\Omega)),\quad u^{\frac{1}{p_N-1}}\in L^{\textbf{p}}_{loc}(0, T; W^{1,\textbf{p}}_{loc}(\Omega)),\\\quad \\
D_i (u^{\frac{1}{p_i-1}} ) \in L^{p_i}_{loc}(0, T; L^{p_i}_{loc}(\Omega)),\,\qquad \text{for all} \quad i=1, ...,N,
\end{cases}
\end{equation}
is a  non-negative, local, weak sub (super)-solution to \eqref{eq1.1}, if for every compact set $E\subset \Omega$ and every sub-interval $[t_1, t_2]\subset (0, T]$
\begin{equation}\label{eq1.3}
\int\limits_E u\,\varphi\,dx \Big|^{t_2}_{t_1}+\int\limits^{t_2}_{t_1}\int\limits_E \Big[u\, \varphi_t+\sum_i A_i(x, t, u, D u)\,D_i \varphi \Big]\,dx dt \leqslant (\geqslant) 0,
\end{equation}
for all non-negative testing functions
$$\varphi\in W^{1,p_N}_{loc}(0, T; L^{p_N}(E))\cap L^{\textbf{p}}_{loc}(0,T; W^{1,\textbf{p}}_0(E)).$$
We say that $u$ is a non-negative, local weak solution to equation \eqref{eq1.1} if $u$ is both a non-negative, local weak
sub- and super-solution.
\noindent

\subsection{Local Energy Estimates}

\begin{lemma}\label{lem2.6}
Let $u$ be a non-negative, local weak sub(super)-solution to  \eqref{eq1.1}-\eqref{eq1.2}. Then there exists a constant $\gamma >0$, depending only on the data, such that, for every cylinder $Q_{r, \eta}(y, \tau)\subset \Omega_T$, every level $k>0$, and every piecewise smooth cutoff function $\zeta$ vanishing on $\partial K_{r}(y)$ and such that $0\leqslant \zeta \leqslant 1$, there holds
\begin{multline}\label{eq2.3}
\sup\limits_{\tau-\eta\leqslant t \leqslant \tau}\int\limits_{K_{r}(y)\times\{t\}}g_{\pm}(u^{\frac{1}{p_N-1}}, k^{\frac{1}{p_N-1}})\,\zeta^{p_N}\,dx\, +\\+\gamma^{-1}  \sum_i\iint\limits_{Q_{r,\eta}(y, \tau)}u^{\frac{p_N-p_i}{p_N-1}}|D_i (u^{\frac{1}{p_N-1}}-k^{\frac{1}{p_N-1}})_{\pm}|^{p_i}\zeta^{p_N}\,dx\,dt\\\leqslant \int\limits_{K_{r}(y)\times\{\tau-\eta\}}g_{\pm}(u^{\frac{1}{p_N-1}}, k^{\frac{1}{p_N-1}})\,\,\zeta^{p_N}\,dx+\gamma\iint\limits_{Q_{r,\eta}(y, \tau)} g_{\pm}(u^{\frac{1}{p_N-1}}, k^{\frac{1}{p_N-1}})\,|\zeta_t|\,dx dt\, +\\+\gamma \sum_i\iint\limits_{Q_{r,\eta}(y, \tau)}u^{\frac{p_N-p_i}{p_N-1}}(u^{\frac{1}{p_N-1}}-k^{\frac{1}{p_N-1}})_{\pm}^{p_i}|D_i \zeta|^{p_i}\,dx\,dt\, .
\end{multline}
\end{lemma} \noindent We postpone the nowadays standard proof of \eqref{eq2.3} to Section \ref{proof-energy-estimates}, in order to leave space to what is really new. At this stage we only remark that the assumption $ u \in L^{\frac{p_N}{p_N-1}}_{loc}(\Omega_T) $ of \eqref{functional} implies that each integral above is convergent. Indeed, as usual, this assumption is crucial for the extrapolation of \eqref{eq2.3} from equation \eqref{eq1.1}-\eqref{eq1.2}. 

\subsection{Local Upper Bound}
Local weak sub-solutions to \eqref{eq1.1}-\eqref{eq1.2} are bounded above quantitatively. This property can be extracted from the study performed in \cite{CiaHenSavSkr1} for De Giorgi classes of functions satisfying  doubly nonlinear parabolic equations that exhibit nonstandard growth conditions. 

\begin{lemma}\label{lem2.7}
Let $u$ be a non-negative, local, weak sub-solution to \eqref{eq1.1}-\eqref{eq1.2}, then for any cylinder $Q_{8\rho, (8\rho)^{p}}(x_0, t_0)\subset \Omega_T$ there holds
\begin{equation}\label{eq2.4}
\sup\limits_{Q_{\rho, \rho^{p}}(x_0, t_0)} u\leqslant \gamma \bigg(\frac{1}{\rho^{N+p}}\iint\limits_{Q_{2\rho,(2\rho)^{p}}(x_0, t_0)} u^{\frac{p_N}{p_N-1}}\, dx dt\bigg)^{\frac{p_N-1}{p_N}},
\end{equation}with a constant $\gamma >0$ depending only on the data.
\end{lemma}

\section{A De Giorgi-type Lemma}\label{Sec.3}
\noindent The following lemma is an anisotropic version of the classic De-Giorgi type lemma (see \cite{DiB, DiBGiaVes1, DiBGiaVes2, DiBGiaVes3}). Although a more general version of it could be extracted from the proof of \cite[Lemma 2.4]{CiaHenSavSkr2}, we decided to report the full proof here, in order to show that for the special case \eqref{eq1.1} there is no need of intrinsic geometry (see \cite{DiB} or \cite{Urb} for more details)  since here the parabolic and the elliptic energy terms suitably balance the De Giorgi type iteration.  

\begin{lemma}\label{lem3.1}
 Let $a \in (0,1)$ and $k>0$ be fixed real numbers, and let $u$ be a non-negative, local weak super-solution to \eqref{eq1.1}-\eqref{eq1.2}. Then there exists $\nu_-\in (0, 1)$, depending only on the data, $a$, $r$ and $\eta$ as in \eqref{eq11.12}, such that if
\begin{equation}\label{eq3.1}
|Q_{r,\eta}(y, \tau)\cap\{u\leqslant k\}|\leqslant \nu_- |Q_{r,\eta}(y, \tau)|,
\end{equation}
then
\begin{equation}\label{eq3.2}
u(x,t)\geqslant a\,k,\quad (x,t)\in Q_{\frac{r}{2}, \frac{\eta}{2}}(y, \tau).
\end{equation}
\end{lemma}
\begin{proof}
For $j=0,1,2, ...$ we define 
$$k_j^{\frac{1}{p_N-1}}=[a\,k]^{\frac{1}{p_N-1}}+(1-a^{\frac{1}{p_N-1}}) k^{\frac{1}{p_N-1}} 2^{-j},$$ 
$$r_j:=\frac{r}{2}\,(1+2^{-j}),\quad \eta_j:=\frac{\eta}{2}(1+2^{-j}),$$
$$K_j:=\K_{r_j}(y),\quad Q_j:=K_j\times I_j,\quad I_j:= (\tau-\eta_j, \tau).$$
Now, we consider piecewise smooth cutoff functions $\zeta_j:=\zeta^{(1)}_j(x) \zeta^{(2)}_j(t)$, where $\zeta^{(1)}_j(x)\in C^1_0(K_j)$ is such that 
$$\zeta^{(1)}_j(x)=1 \quad \text{for} \quad  x\in K_{j+1}, \qquad  0\leqslant \zeta^{(1)}_j(x)\leqslant 1\,,$$
$$\text{and}\qquad |D_i \zeta^{(1)}_j(x)|^{p_i}\leqslant \frac{\gamma 2^{j}}{r^p}, i=1, ..., N,$$
while $\zeta^{(2)}_j(t)\in C^1(\mathbb{R})$ is such that
$$\zeta^{(2)}_j(t)=0 \quad \text{for} \quad  t\leqslant \tau-\eta_j; \qquad   
\zeta^{(2)}_j(t)=1, \quad \text{for} \quad 
 t\geqslant \tau-\eta_{j+1}, $$
 $$0\leqslant \zeta^{(2)}_j(t)\leqslant 1, \quad \text{verifying} \quad  |\zeta^{(2)}_{j, t}|\leqslant \frac{\gamma 2^j}{\eta}.$$
Using Lemmas \ref{lem2.5} and \ref{lem2.6}, we obtain
\begin{multline*}
\sup\limits_{t \in I_j}\int\limits_{K_j \times\{t\}}g_-(u^{\frac{1}{p_N-1}}, k^{\frac{1}{p_N-1}}_j)\,\zeta^{p_{N}}_j\,dx+\gamma^{-1}\sum_i \iint\limits_{Q_j}u^{\frac{p_N-p_i}{p_N-1}}|D_i(u^{\frac{1}{p_N-1}}-k^{\frac{1}{p_N-1}}_j)_{-}|^{p_{i}}\,\zeta^{p_{N}}_j dx\,dt\\\leqslant
\frac{\gamma 2^{j\gamma}}{\eta}\iint\limits_{Q_j}g_-(u^{\frac{1}{p_N-1}}, k^{\frac{1}{p_N-1}}_j)\,dx dt+\frac{\gamma 2^{j \gamma}}{r^p}\sum_i \iint\limits_{Q_j}u^{\frac{p_N-p_i}{p_N-1}}(u^{\frac{1}{p_N-1}}-k^{\frac{1}{p_N-1}}_j)_{-}^{p_{i}}\,dx\,dt \\ \leqslant \gamma 2^{j \gamma}\,k^{\frac{p_N}{p_N-1}}\Big(\frac{1}{\eta}+\frac{1}{r^p}\Big)|Q_j\cap\{u\leqslant k_j\}|.
\end{multline*}
To estimate the left-hand side of this inequality we define
$$v:=\max(u, a\,k).$$
and note that each one of the terms on the left-hand side can be estimated in term of $v$ as follows:
\begin{multline*}
\int\limits_{K_j\times\{t\}}g_-(u^{\frac{1}{p_N-1}}, k^{\frac{1}{p_N-1}}_j)\,\zeta^{p_{N}}_j\,dx\geqslant
\int\limits_{K_j\times\{t\}}g_-(v^{\frac{1}{p_N-1}}, k^{\frac{1}{p_N-1}}_j)\,\zeta^{p_{N}}_j\,dx\geqslant\\\geqslant  \frac{1}{2}(ak)^{\frac{p_N-2}{p_N-1}}
\int\limits_{K_j\times\{t\}}(v^{\frac{1}{p_N-1}}-k^{\frac{1}{p_N-1}}_j)^2_-\,\zeta^{p_N}_j\,dx,
\end{multline*}
and also, for $\alpha_i={p_N}/{p_i}$,
\begin{multline*}
\iint\limits_{Q_j} u^{\frac{p_N-p_i}{p_N-1}}|D_i(u^{\frac{1}{p_N-1}}-k^{\frac{1}{p_N-1}}_j)_{-}|^{p_i}\,\zeta_j^{p_{N}}\,dxdt\\\geqslant \iint\limits_{Q_j\cap\{u\geqslant ak\}}u^{\frac{p_N-p_i}{p_N-1}}|D_i(u^{\frac{1}{p_N-1}}-k^{\frac{1}{p_N-1}}_j)_{-}|^{p_i}\,\zeta_j^{p_{N}}\,dxdt \\ =
\iint\limits_{Q_j\cap\{u\geqslant ak\}}v^{\frac{p_N-p_i}{p_N-1}}|D_i(u^{\frac{1}{p_N-1}}-k^{\frac{1}{p_N-1}}_j)_{-}|^{p_i}\,\zeta_j^{p_{N}}\,dxdt\\ \geqslant
\gamma^{-1}a^{\frac{p_N-p_i}{p_N-1}}\iint\limits_{Q_j}|D_i(v^{\frac{1}{p_N-1}}-k^{\frac{1}{p_N-1}}_j)^{\alpha_i}_{-}|^{p_i}\,\zeta_j^{p_{N}}\,dxdt.
\end{multline*}    

Hence, collecting the previous estimates we arrive at
\begin{multline*}
(ak)^{\frac{p_N-2}{p_N-1}}\,\sup\limits_{t\in I_j}
\int\limits_{K_j\times\{t\}}(v^{\frac{1}{p_N-1}}-k^{\frac{1}{p_N-1}}_j)^2_-\,\zeta^{p_N}_j\,dx
+\gamma^{-1}a^{\frac{p_N-p_{1}}{p_N-1}}
\sum_i\iint\limits_{Q_j}|D_i(v^{\frac{1}{p_N-1}}-k^{\frac{1}{p_N-1}}_j)^{\alpha_i}_{-}|^{p_i}\,\zeta_j^{p_{N}}\,dxdt
\\
\leqslant\gamma 2^{j \gamma}\,k^{\frac{p_N}{p_N-1}}\Big(\frac{1}{\eta}+\frac{1}{r^p}\Big)|Q_j\cap\{u\leqslant k_j\}|. 
\end{multline*}
By means of H\"{o}lder's inequality together with Lemma \ref{lem2.2}, for $\beta=2$  and $\alpha=\frac{1}{N}\sum\limits_{i=1}^N\alpha_i=\frac{p_N}{p}$, we obtain
\[
\left(k^{\frac{1}{p_N-1}}_j-k^{\frac{1}{p_N-1}}_{j+1}\right)^{p_N} |Q_{j+1}\cap\{u\leqslant k_{j+1}\}|=\left(k^{\frac{1}{p_N-1}}_j-k^{\frac{1}{p_N-1}}_{j+1}\right)^{\alpha p} |Q_{j+1}\cap\{v\leqslant k_{j+1}\}| \]
\begin{eqnarray*}
& \leqslant &
\iint\limits_{Q_j}\left[(v^{\frac{1}{p_N-1}}-k^{\frac{1}{p_N-1}}_j)_{-}\zeta^{p_{N}}_j \right]^{\alpha p}\,dxdt\\
& \leqslant & \Big(\iint\limits_{Q_j}\big[(v^{\frac{1}{p_N-1}}-k^{\frac{1}{p_N-1}}_j)_{-}\zeta^{p_{N}}_j\big]^{p\frac{\alpha N+2}{N}}\,dxdt\Big)^{\frac{N}{N+2/\alpha}}|Q_j\cap\{v\leqslant k_j\}|^{\frac{2/\alpha}{N+2/\alpha}}\\ 
& \leqslant & 
\gamma\Big(\sup\limits_{t\in I_j}\int\limits_{K_j\times\{t\}}(v^{\frac{1}{p_N-1}}-k^{\frac{1}{p_N-1}})^{2}_{-}\zeta^{p_{N}}_j\,dx\Big)^{\frac{p}{N+2/\alpha}}\times\\
& & \times \Big(\sum_i \iint\limits_{Q_j}|D_i \big[(v^{\frac{1}{p_N-1}}-k^{\frac{1}{p_N-1}}_j)^{\alpha_i}_{-}\zeta^{p_{N}}_j\big]|^{p_i}\,dxdt\Big)^{\frac{N}{N+2/\alpha}}|Q_j\cap\{u\leqslant k_j\}|^{\frac{2/\alpha}{N+2/\alpha}}\\ 
& \leqslant &\gamma 2^{j\gamma} a^{-\gamma}
\,k^{\frac{p_N}{p_N-1}}\Big(\frac{1}{\eta}+\frac{1}{r^p}\Big)^{\frac{N+p}{N+2/\alpha}}|Q_j\cap\{u\leqslant k_j\}|^{1+\frac{p}{N+2/\alpha}}.
\end{eqnarray*}
Set $y_j:=\dfrac{|Q_j\cap\{u\leqslant k_j\}|}{|Q_j|}$, then from the previous we obtain
\begin{equation*}
y_{j+1}\leqslant \gamma 2^{j\gamma} a^{-\gamma}\Big(\frac{r^p}{\eta}\Big)^{\frac{N}{N+2/\alpha}}\Big(1+ \frac{\eta}{r^p}\Big)^{\frac{N+p}{N+2/\alpha}}y_{j}^{1+\frac{p}{N+2/\alpha}},
\end{equation*}
which yields $\lim\limits_{j\rightarrow \infty}y_j=0$, provided that $y_0\leqslant \nu_-$, where $\nu_-$ is chosen to satisfy
\begin{equation}\label{eq11.12}
\nu_-:=\frac{a^{\gamma}}{\gamma}\Big(\frac{\eta}{r^p}\Big)^{\frac{N}{p}}\Big(1+ \frac{\eta}{r^p}\Big)^{-\frac{N+p}{p}},
\end{equation}
therefore proving proves \eqref{eq3.2}.
\end{proof}

\section{Expansion of Positivity}\label{Sec.4}

\noindent Expansion of positivity is the main and crucial result to obtain Harnack's inequality, and reads as follows.

\begin{propo}\label{pr4.1}
Let $u$ be a non-negative, local weak super-solution to \eqref{eq1.1}-\eqref{eq1.2} and let us assume that for some $r$, $k>0$ and $\alpha_0 \in (0,1)$ the following measure information is at stake
\begin{equation}\label{eq4.1}
|\K_r(y)\cap\{u(\cdot, s)\geqslant k\}|\geqslant \alpha_0 |\K_r(y)|.
\end{equation}
Then, there exist numbers $\eta_0$, $\delta_0 \in (0, 1)$ depending only on the data and $\alpha_0$ such that
\begin{equation}\label{eq4.2}
u(x, t)\geqslant \eta_0\,k,\quad x\in \K_{2r}(y),
\end{equation}
for all times
\begin{equation}\label{eq4.3}
s+\frac{\delta_0}{2}r^p\leqslant t \leqslant s+\delta_0 r^p,
\end{equation}
provided that 
$$\K_{8r}(y)\times(s, s+r^p)\subset \Omega_T.$$
\end{propo}
\noindent By repeated application of Proposition \ref{pr4.1} we obtain
\begin{corollary}\label{cor4.1}
Under conditions stated in Proposition \ref{pr4.1}, for any $A>0$ there exists number $\bar{\eta}\in (0, 1)$ depending only on the data, $\alpha_0$  and $A$ such that 
\begin{equation*}
u(x, t)\geqslant \bar{\eta}\,k,\quad x\in \K_{2r}(y),
\end{equation*}
for all $t$ 
\begin{equation*}
s+\frac{A}{2}r^p\leqslant t \leqslant s+A r^p,
\end{equation*}
provided that \eqref{eq4.1} holds and $\K_{8 r}(y)\times (s, s+ A(8r)^{p}) \subset \Omega_T.$
\end{corollary}


\noindent The proof of Proposition \ref{pr4.1} is performed in two major steps, given by Lemma \ref{lem4.1} and Lemma \ref{lem4.2}.

\noindent We start by introducing the change of variables $$t \rightarrow \frac{t-s}{r^{p}}= \tau,\quad x_i \rightarrow \frac{x_i-y_i}{r^{\frac{p}{p_i}}}=z_i,\quad i=1, ..., N,$$
which maps the cylinder $\K_{8r}\times(s, s+r^{p})$ into $Q:=K_8(0)\times (0,1)$. Function $u$ keeps on being a local weak super-solution 
to \eqref{eq1.1}-\eqref{eq1.2}, in $Q$, and inequality \eqref{eq4.1} translates into
\begin{equation}\label{eq4.4}
|K_1(0)\cap \{u(\cdot, 0) \geqslant k\}|\geqslant \alpha_0 |K_1(0)|.
\end{equation}

\subsection{Propagation of Positivity in time}

\begin{lemma}\label{lem4.1}
Let $u$ be a non-negative, local weak super-solution to \eqref{eq1.1}-\eqref{eq1.2} in $Q$ and assume that 
inequality \eqref{eq4.4} holds. 
 Then there exist $\varepsilon_1$, $\delta_1 \in (0, 1)$, depending only on the data and on $\alpha_0$, such that
\begin{equation}\label{eq4.5}
|K_1(0)\cap\{u(\cdot, \tau)\geqslant \varepsilon_1\,k\}|\geqslant \frac{\alpha_0}{4}|K_1(0)|,
\end{equation}
for all 
\begin{equation*}
0 \leqslant \tau\leqslant \delta_1.
\end{equation*}
\end{lemma}

\begin{proof}
The proof relies on the energy estimates of Lemma \ref{lem2.6}, written over the cylinder $K_1(0) \times (0, \delta_1]$, where the cutoff function $\zeta\in C^1_0(K_1(0))$ can be chosen as a time-independent function satisfying for a fixed $\sigma \in (0,1)$ the usual properties
\[\zeta(x)=1 \, \, \text{for} \, \, x \in K_{1-\sigma}(0), \qquad 0\leqslant \zeta(x)\leqslant 1, \qquad |D_i \zeta(x)|^{p_i}\leqslant \frac{\gamma}{\sigma^{p_N}} \qquad  \forall i=1, ...,N\,.\] Now, by observing that, for any $0\leqslant u<k$ we can reduce from above the following term
\begin{equation*}
    \begin{aligned}
g_-(u^{\frac{1}{p_N-1}}, k^{\frac{1}{p_N-1}})&=(p_N-1)\int\limits^{k^{\frac{1}{p_N-1}}}_{u^{\frac{1}{p_N-1}}} s^{p_N-2}(k^{\frac{1}{p_N-1}}-s) ds=\\
&=\frac{1}{p_N} k^{\frac{p_N}{p_N-1}}-u \bigg(k^{\frac{1}{p_N-1}}-\frac{p_N-1}{p_N}u^{\frac{1}{p_N-1}}\bigg)\leqslant  \frac{k^{\frac{p_N}{p_N-1}}}{p_N},
    \end{aligned}
\end{equation*}
and when restricting on the set $\{u \leqslant \varepsilon_1 k\}$ for a small $0< \varepsilon_1 < 1/p_N$, we can reduce it from below 
\[ g_-(u^{\frac{1}{p_N-1}}, k^{\frac{1}{p_N-1}})\geqslant \frac{1}{p_N} k^{\frac{p_N}{p_N-1}}-\varepsilon_1 k^{\frac{p_N}{p_N-1}}=k^{\frac{p_N}{p_N-1}} (1/p_N-\varepsilon_1 )\,,\]
we arrive at the following estimate, valid for all times $0\leqslant  \tau \leqslant \delta_1$:
\begin{equation*}
    \begin{aligned}
\, k^{\frac{p_N}{p_N-1}} & (1/p_N-\varepsilon_1) |{K}_{1-\sigma}(0)\cap\{u(\cdot, \tau)\leqslant \varepsilon_1 k\}|\leqslant
\sup\limits_{0\leqslant \tau\leqslant \delta_1}\int\limits_{K_1(0)}g_-(u^{\frac{1}{p_N-1}}, k^{\frac{1}{p_N-1}})\zeta^{p_N}\,dz \\
& \leqslant \int\limits_{K_1(0)\times\{0\}}g_-(u^{\frac{1}{p_N-1}}, k^{\frac{1}{p_N-1}})\,dz+\frac{\gamma}{\sigma^{p_N}}\sum_i\int\limits_{0}^{\delta_1 }\int\limits_{K_1(0)}u^{\frac{p_N-p_i}{p_N-1}}(u^{\frac{1}{p_N-1}}-k^{\frac{1}{p_N-1}})^{p_i}_-\,dz d\tau  \\
&\qquad \leqslant  \bigg[  \frac{1-\alpha_0}{p_N}+\frac{\gamma \delta_1}{\sigma^{p_N}}\bigg]k^{\frac{p_N}{p_N-1}}|K_1(0)|\, ,
    \end{aligned}
\end{equation*}
which by the measure estimate 
\[|K_{1-\sigma}(0)\cap\{u(\cdot, \tau)\leqslant \varepsilon_1 k\}|\geqslant |K_1(0)\cap\{u(\cdot, t)\leqslant \varepsilon_1 k\}|-N \sigma |K_1(0)|,\]
implies an estimate from above of the sub-level sets of $u$ as
\begin{equation*}
|K_1(0)\cap\{u(\cdot, \tau)\leqslant \varepsilon_1 k\}|\leqslant \bigg[N\,\sigma +\frac{p_N}{1-\varepsilon_1 p_N}\Big( \frac{1-\alpha_0}{p_N}+\frac{\gamma\,\delta_1}{\sigma^{p_N}}\Big)\bigg]\,|K_1(0)|.
\end{equation*}
Now, by choosing
\begin{equation*}
\sigma=\frac{\alpha_0}{8 N},\qquad  \varepsilon_1= \frac{\alpha_0 }{(2-\alpha_0)p_N},\qquad \delta_1=\frac{\alpha_0(1-\alpha_0)\sigma^{p_N}}{16\gamma p_N},
\end{equation*}
we arrive at the required \eqref{eq4.5}, which proves the lemma.
\end{proof}

\subsection{A Measure Shrinking Lemma}
\begin{lemma}\label{lem4.2}
Assume that the conditions of Lemma \ref{lem4.1} are  fulfilled and let $\delta_1, \varepsilon_1, \alpha_0$ given in lemma \ref{lem4.1}. Then, for any $\nu \in (0, 1)$ there exists $j_*$, depending only on the data, $\varepsilon_1$, $\alpha_0$ and $\nu$, such that
\begin{equation}\label{eq4.7}
\left|K_4(0)\times\left(\frac{\delta_1}{4}, \delta_1 \right)\cap\{u\leqslant \frac{k}{2^{j_*}}\}\right|\leqslant  \frac{3}{4}\nu\,\delta_1 |K_4(0)|.
\end{equation}
\end{lemma}
\begin{proof}Let $j_*\in \mathbb{N}$ to be determined and, for $j \in \mathbb{N}$ such that $1\leqslant j\leqslant j_*$, let us set $k_j:=\varepsilon_1^{j} k$ and use Lemma \ref{lem2.6} for the pair of cylinders 
\[Q_1:=K_{4}(0)\times \bigg(\frac{\delta_1}{4},\, \delta_1 \bigg) \qquad \text{and} \qquad Q_2:=K_{6}(0)\times\bigg(0, \, \delta_1\bigg)\,,\] with a piecewise smooth cutoff function $\zeta(x, t)$ vanishing on the parabolic boundary of $Q_2$ and verifying \[|\zeta_{x_i}|^{p_i}\leqslant \gamma\,, \quad \forall i=1, ..., N, \qquad \text{and} \qquad |\zeta_t|\leqslant \frac{\gamma}{\delta_1}\,,\] to obtain
\[
\sum_i \iint\limits_{Q_1\cap\{u\leqslant k_j\}}u^{\frac{p_N}{p_N-1}-p_i}|D_i u|^{p_i}\,dz\,d\tau\leqslant\gamma
\sum_i \iint\limits_{Q_2}u^{\frac{p_N-p_i}{p_N-1}}|D_i(u^{\frac{1}{p_N-1}}-k^{\frac{1}{p_N-1}}_j)_{-}|^{p_{i}}\zeta^{p_{N}}\,dz\,d\tau \] 
\begin{eqnarray*}
& \leqslant & \gamma\iint\limits_{Q_2}\,g_-(u^{\frac{1}{p_N-1}}, k^{\frac{1}{p_N-1}}_j)\,|\zeta_t|\,dz d\tau+\gamma \sum_i\iint\limits_{Q_2}u^{\frac{p_N-p_i}{p_N-1}}(u^{\frac{1}{p_N-1}}-k^{\frac{1}{p_N-1}}_j)_{-}^{p_{i}}|D_i\zeta|^{p_{i}}\,dz\,d\tau
\\
& \leqslant &\frac{\gamma}{\delta_1} \,k^{\frac{p_N}{p_N-1}}_j |Q_1|.
\end{eqnarray*}
Thanks to the discrete iso-perimetric inequality Lemma \ref{lem2.3} and the assumption \eqref{eq4.4} we arrive at
\[
k_{j+1}\,|K_{4}(0)\cap \{u(\cdot, \tau)\leqslant k_{j+1}\}|\leqslant \frac{\gamma}{\alpha_0}\sum_i\,\int\limits_{K_{4}(0)\cap\{k_{j+1}\leqslant u(\cdot, \tau)\leqslant k_j\}}|D_i u|\,dz
\]
for every time  
$\tau\in \big(\frac{\delta_1}{4}, \delta_1\big)$. Integrating this inequality over $\big(\frac{\delta_1}{4}, \delta_1\big)$ by usual manipulations we get
\begin{multline*}
\varepsilon_1 k_{j}\,|Q_1\cap\{u\leqslant k_{j+1}\}|\leqslant \frac{\gamma}{\alpha_0}\sum_i\,\iint\limits_{Q_1\cap\{k_{j+1}\leqslant u\leqslant k_j\}}|D_i u|\,dz d\tau\\\leqslant \frac{\gamma}{\alpha_0}\sum_i\,\Big(\iint\limits_{Q_1\cap\{ u\leqslant k_j\}}u^{\frac{p_N}{p_N-1}-p_i}|D_i u|^{p_i}\,dz d\tau\Big)^{\frac{1}{p_i}}\Big(\iint\limits_{Q_1\cap\{k_{j+1}\leqslant u\leqslant k_j\}}u^{\frac{p_i-\frac{p_N}{p_N-1}}{p_i-1}}\,dx dt\Big)^{1-\frac{1}{p_i}}\\\leqslant \frac{\gamma}{\alpha_0} \sum_i\,\Big(k^{\frac{p_N}{p_N-1}}_j |Q_1|/\delta_1\Big)^{\frac{1}{p_i}}\Big(k^{\frac{p_i-\frac{p_N}{p_N-1}}{p_i-1}}_j|Q_1\cap\{k_{j+1}\leqslant u\leqslant k_j\}|\Big)^{1-\frac{1}{p_i}}\\=\frac{\gamma}{\alpha_0} k_j \sum_i|K_4|^{\frac{1}{p_i}}\,|Q_1\cap\{k_{j+1}\leqslant u\leqslant k_j\}|^{1-\frac{1}{p_i}}\\\leqslant
\gamma \alpha_0^{-1}\, k_j |Q_1|\,\Big(\frac{|Q_1\cap\{k_{j+1}\leqslant u\leqslant k_j\}|}{|Q_1|}\Big)^{1-\frac{1}{p_{1}}}.
\end{multline*}
Thereby we end up with an estimate of the sub-level sets of $u$ such as
\begin{equation*}
\Big(\frac{|Q_1\cap\{u\leqslant k_{j_*}\}|}{|Q_1|}\Big)^{\frac{p_{1}}{p_{1}-1}}\leqslant \gamma \alpha_0^{-2\frac{p_1}{p_1-1}}\,\frac{|Q_1\cap\{k_{j+1}\leqslant u\leqslant k_j\}|}{|Q_1|}.
\end{equation*}
Summing up the last inequalities over $j$, for $1\leqslant j \leqslant j_*$, and using the properties of Lebesgue measure we deduce the shrinking property
\begin{equation*}
|Q_1\cap\{u\leqslant k_{j_*}\}|\leqslant \frac{\gamma \alpha_0^{-2}}{(j_*-1)^{1-\frac{1}{p_{1}}}}|Q_1|.
\end{equation*}
The proof is completed once we choose $j_*$ large enough such that
\[\frac{\gamma\,  \alpha_0^{-2}}{(j_*-1)^{1-\frac{1}{p_{1}}}}\leq \nu \,, \quad \quad  \text{and} \quad \quad j_*= \bigg{\lceil} 1+  \bigg(\frac{\gamma}{\, \nu\,  \alpha_0^{2}} \bigg)^{\frac{p_1}{p_1-1}} \bigg{\rceil} \,,\]
where in the last definition we used the upper integer part. 
\end{proof}

\noindent De Giorgi type Lemma \ref{lem3.1}, for $a=\frac{1}{2}$, and Lemma \ref{lem4.2} yield
\begin{equation*}
u(z, \tau)\geqslant k \varepsilon_1^{j_*}/2,\quad (z, \tau)\in K_2(0)\times \bigg(\frac{\delta_1}{2}, \, \delta_1 \bigg).
\end{equation*}
Returning to the original variables, we arrive at the required \eqref{eq4.2}-\eqref{eq4.3} with $\delta_0=\delta_1$ and $\eta_0=\varepsilon_1^{j_*}/{2}$, completing the proof of Proposition \ref{pr4.1}.

\section{Harnack's Inequality}\label{Sec.5}

\noindent Let us fix $(x_0, t_0)\in \Omega_T$ for which $u(x_0, t_0)>0$ and construct the cylinder \[\K_{8\rho}(x_0)\times (t_0-\bar{C}(8\rho)^{p}, t_0+\bar{C}(8\rho)^{p})\subset \Omega_T\]  with $\bar{C}>0$ to be chosen.
The change of variables
$$t\rightarrow \frac{t-t_0}{\rho^{p}}= \tilde{t},\quad x_i\rightarrow \frac{x_i-x_{0,i}}{\rho^{\frac{p}{p_i}}}=z_i,\quad i=1, ...,N$$
maps this cylinder into $Q:=K_8(0)\times (-\bar{C}8^{p}, \bar{C}8^{p})$, and  turns $v(z,\tilde{t})=\frac{1}{u(x_0, t_0)}\,u (x,y)$ in a bounded, non-negative, weak solution to \eqref{eq1.1}-\eqref{eq1.2} in $Q$, verifying $v(0,0)=1$. 

\vskip0.2cm

\noindent Consider the cylinder $Q_{1,1}(0, 0)$ and fix $\beta>0$, to be determined later. By Lemma \ref{lem2.4}, there exist $\gamma_0=\gamma_0(\beta)>1$,  $(y, \tau)\in Q_{1,1}(0, 0)$  and $r > 0$ such that
\begin{equation}\label{eq5.1}
Q_{r, r^p}(y, \tau)\subset Q_{1,1}(0, 0),\quad r^\beta \sup\limits_{Q_{r, r^p}(y, \tau)}v \leqslant \gamma_0,\quad
r^\beta v(y, \tau)\geqslant \frac{1}{\gamma_0}.
\end{equation}
\noindent \textbf{Claim $1$.} There exist $\epsilon$, $\alpha \in (0, 1)$ depending only on the data and $\beta$ such that
\begin{equation}\label{eq5.2}
\big|Q_{r, r^{p}}(y, \tau)\cap \big\{v\geqslant \epsilon r^{-\beta}\big\}\big|\geqslant \alpha |Q_{r, r^{p}}(y, \tau)|.
\end{equation}
\begin{proof}
Since $v$ is a a bounded, non-negative, weak solution to \eqref{eq1.1}-\eqref{eq1.2}, from Lemma \ref{lem2.7} together with  \eqref{eq5.1} and the fact that
$|Q_{r, r^{p}}(y, \tau)|=r^{N+p}$
we get
\begin{multline*}
\frac{1}{\gamma_0}r^{-\beta} \leqslant v(y, \tau)\leqslant \gamma \bigg(\frac{1}{r^{N+p}}\iint\limits_{Q_{r, r^p}(y, \tau)} v^{\frac{p_N}{p_N-1}}\,dx\,dt\bigg)^{\frac{p_N-1}{p_N}}\leqslant \gamma \epsilon r^{-\beta}+\\+\gamma \bigg(\frac{1}{r^{N+p}}\iint\limits_{Q_{r, r^{p}}(y, \tau)\cap \{v\geqslant \epsilon r^{-\beta}\}} v^{\frac{p_N}{p_N-1}}\,dx\,dt\bigg)^{\frac{p_N-1}{p_N}}\\\leqslant \gamma \epsilon r^{-\beta}+\gamma \gamma_0 r^{-\beta}\bigg(\frac{|Q_{r, r^{p}}(y, \tau)\cap \{v\geqslant \epsilon r^{-\beta}\}|}{|Q_{r, r^{p}}(y, \tau)|}\bigg)^{\frac{p_N-1}{p_N}}.
\end{multline*}
We arrive at the required \eqref{eq5.2}, with
$\alpha=\Big(\dfrac{1}{2 \gamma \gamma^2_0}\Big)^{\frac{p_N}{p_N-1}}$, once we choose  $\epsilon=\dfrac{1}{2\gamma \gamma_0}$. 
\end{proof}

\noindent Inequality \eqref{eq5.2} yields
\begin{equation}\label{eq5.3}
\big|K_{r}(y)\cap\big\{v(\cdot, s)\geqslant \epsilon r^{-\beta}\big\}\big|\geqslant \alpha |K_r(y)|,\quad\text{for some}\quad s\in (\tau -r^{p}, \tau),
\end{equation}
and then, by applying Proposition \ref{pr4.1} (with $k=\epsilon r^{-\beta}$ and $\alpha_0=\alpha$), we get
\begin{equation}\label{eq5.4}
v(z, t_0)\geqslant \eta_0\,\epsilon r^{-\beta},\quad x\in K_{2r}(y),\quad t_0=s+\delta_0\,r^{p},
\end{equation}
for some $\eta_0$, $\delta_0\in (0, 1)$ depending only on the data and on $\beta$.

\vskip0.2cm 

\noindent We can now repeat this procedure (one and again) starting from \eqref{eq5.4}. In fact, the expansion of positivity now applied for level
$\eta_0\,\epsilon r^{-\beta}$, radius  $2r$ and $\alpha_0=1$ yields 
\begin{equation}\label{eq5.5}
v(z, t_1)\geqslant \bar{\eta}\,\eta_0\epsilon r^{-\beta},\quad x \in K_{2^2r}(y),\quad t_1:= t_0+\bar{\delta}(2r)^{p},
\end{equation}
where $\bar{\eta}$, $\bar{\delta}\in (0, 1)$ depend only on the data and are independent of $r$ and $\beta$. By an iterative scheme one arrives at 
\begin{equation}\label{eq5.6}
v(z,t_n)\geqslant \bar{\eta}^n\, \eta_0\epsilon r^{-\beta},\quad x \in K_{2^{n+1} r}(y),\quad t_n:= t_{n-1}+\bar{\delta}(2^n r)^{p}.
\end{equation}
Choosing $n$ so large that $2^{n+1}\,r \geqslant 2$, the cube $K_{2^{n+1} r}(y)$ covers the cube $K_1(0)$ and then
\begin{equation}\label{eq5.7}
v(z,t_n)\geqslant \bar{\eta}^n\, \eta_0\epsilon r^{-\beta},\quad x \in K_{1}(0)
\end{equation}
where time level $t_n$ is given by
$$t_n=t_0+\bar{\delta}\sum\limits_{j=1}^n (2^j r)^{p}=s+\delta_0 r^{p}+\bar{\delta}\sum\limits_{j=1}^n (2^j r)^{p},$$
which can be estimated as
$$r^{p}\bar{\delta} 2^{n p}-1\leqslant t_n \leqslant 1+ \bar{\delta}(2^{n+1} r)^{p}.$$
Choosing $n$ such that $2\leqslant \bar{\delta} (2^n r)^{p}\leqslant 2^{p+1}$, from \eqref{eq5.7} we obtain
\begin{equation}\label{eq5.8}
v(z, \bar{t}) \geqslant \bar{\eta}^n\, \eta_0\epsilon r^{-\beta} \geqslant \eta_0\epsilon \Big(\frac{\bar{\delta}}{2^{p+1}}\Big)^{\frac{\beta}{p}}\big(2^\beta\,\bar{\eta}\big)^n = \eta_0\, \epsilon\,  \Big(\frac{\bar{\delta}}{2^{p+1}}\Big)^{\frac{-\log_2(\bar{\eta})}{p}} \ , \quad x \in K_{1}(0), \quad  1 \leqslant \bar{t}\leqslant 2^{p+2} \ , 
\end{equation}
once we take $2^\beta\,\bar{\eta}=1$. Returning to the original variables shows inequality \eqref{eq1.5} forward in time.

{\subsection{Backward Estimate}
Let $(x_0, t_0)\in \Omega_T$ be such that $\K_{8\rho}(x_0)\times (t_0-\bar{C}(8\rho)^{p}, t_0+\bar{C}(8\rho)^{p})\subset \Omega_T$, and let $y$ be an arbitrary point in $\mathcal{K}_\rho(x_0)$.  
Applying the forward Harnack inequality in \eqref{eq1.5}, which has already been established, with $(x_0,t_0)$ replaced by $(y,t_0-\bar{C}\rho^p)$, we obtain
\begin{equation}\label{BE_1}
u(y, t_0-\bar{C}\rho^p)\leqslant C\,\inf\limits_{\mathcal{K}_\rho(y)}u(\cdot, t_0)\leqslant C\,u(x_0, t_0).
\end{equation}
By definition of supremum, for any $\epsilon \in (0,1)$, there exists $y_\epsilon\in \mathcal{K}_\rho(x_0)$ such that
\[
\sup\limits_{\mathcal{K}_\rho(x_0)} u(\cdot, t_0-\bar{C}\rho^p)
\leqslant u(y_\epsilon, t_0-\bar{C}\rho^p)+\epsilon.
\]
Since $y$ is arbitrary in $\mathcal{K}_\rho(x_0)$, estimate \eqref{BE_1} applies in particular to $y_\epsilon$, which yields
\[
\sup\limits_{\mathcal{K}_\rho(x_0)} u(\cdot, t_0-\bar{C}\rho^p)
\leqslant C\,u(x_0, t_0)+\epsilon.
\]
Letting $\epsilon \to 0$ shows the backward Harnack inequality in \eqref{eq1.5}, completing the proof of Theorem~\ref{th1.1}.

\section{H\"{o}lder Continuity}\label{Sec.6}
\noindent In this section we prove Theorem \ref{th1.2}. let us fix a point $(x_0, t_0)\in \Omega_T$ and let $R>0$ be such that $Q_{8R,(8R)^{p}}(x_0, t_0)\subset \Omega_T.$
If for any $0<\rho < R$ 
$$\osc\limits_{Q_{\rho,  \rho^{p}}(x_0, t_0)} u\leqslant \rho, $$
then there is nothing to prove. So, assume that there
exists $0<\rho_0< R$ such that 
$$\osc\limits_{Q_{\rho_0,  \rho^{p}_0}(x_0, t_0)} u\geqslant \rho_0 $$
and set 
$$Q_{\rho_0}:=Q_{\rho_0, \rho^{p}_0}(x_0, t_0),\qquad \mu^+_0:=\sup\limits_{Q_{\rho_0}}u,\qquad \mu^-_0:=\inf\limits_{Q_{\rho_0}}u,\qquad \omega_0:=\mu^+_0-\mu^-_0.$$
We will construct the decreasing sequences of numbers  \[\rho_j:=\epsilon^{j}\,\rho_0,  \qquad \omega_j:=\sigma^j \omega_0 \, \qquad \mbox{for} \ \ \epsilon, \sigma\in (0, 1) \]
and cylinders $ Q_{\rho_j}:=Q_{\rho_j,\rho^{p}_j}(x_0, t_0)$ such that
\begin{equation}\label{eq6.1}
\mu^+_j:=\sup\limits_{Q_{\rho_j}}u,\qquad \mu^-_j:=\inf\limits_{Q_{\rho_j}} u,\qquad \mu^+_j-\mu^-_j\leqslant \omega_j,\qquad j=0, 1, 2, ....
\end{equation}
Assume first that
\begin{equation}\label{eq6.2}
\mu^+_0\leqslant \big(1+\eta\big)\omega_0,
\end{equation}
where $\eta\in (0, 1)$ is to be determined depending only on the data. Observe that one of the following cases must be true : either
\begin{equation}\label{eq6.3}
\big|Q_{\frac{3}{4}\rho_0}\cap\big\{u\geqslant \mu^+_0-\frac{1}{2}\omega_0\big\}\big|=0,
\end{equation}
or
\begin{equation}\label{eq6.4}
\big|Q_{\frac{3}{4}\rho_0}\cap\big\{u\geqslant \mu^+_0-\frac{1}{2}\omega_0\big\}\big|>0.
\end{equation}
In the first case, we immediately get 
$$u(x,t)\leqslant \mu^+_0-\frac{1}{2}\omega_0,\qquad \forall (x,t)\in Q_{\frac{3}{4}\rho_0},$$
obtaining therefore the reduction of the oscillation
\begin{equation}\label{eq6.5}
\osc\limits_{Q_{\frac{1}{2}\rho_0}}u\leqslant \frac{1}{2}\,\omega_0.
\end{equation}
Assume now that \eqref{eq6.4} is in force.
\vskip0.2cm \noindent \textbf{Claim $1$.} There exist $\xi$, $\alpha\in (0, 1)$ depending only on the data  such that
\begin{equation}\label{eq6.6}
\big|Q_{\rho_0}\cap\big\{u\geqslant \xi\, \omega_0\big\}\big|\geqslant \alpha |Q_{\rho_0}|.
\end{equation}
\begin{proof}
From \eqref{eq6.4} and Lemma \ref{lem2.7}, applied to the pair of cylinders $Q_{\frac{3}{4}\rho_0}$ and $Q_{\rho_0}$, and recalling \eqref{eq6.2}, we obtain 
\begin{eqnarray*}
\frac{1}{4}\omega_0 & \leqslant & \mu_0^+ - \frac{3}{4} \omega_0 \leqslant \sup\limits_{Q_{\frac{3}{4}\rho_0}}u\leqslant \gamma\,\bigg(\frac{1}{\rho^{N+p}_0}\iint\limits_{Q_{\rho_0}} u^{\frac{p_N}{p_N-1}}\,dx dt\bigg)^{\frac{p_N-1}{p_N}} \\
& \leqslant & \gamma \,\xi\,\omega_0
+ \gamma\,\bigg(\frac{1}{\rho^{N+p}_0}\iint\limits_{Q_{\rho_0}\cap\{u\geqslant \xi\,\omega_0\}} u^{\frac{p_N}{p_N-1}}\,dx dt\bigg)^{\frac{p_N-1}{p_N}} \\
& \leqslant & \gamma \xi \omega_0+\gamma \mu^+_0\bigg(\frac{|Q_{\rho_0}\cap\{u\geqslant \xi \omega_0\}|}{|Q_{\rho_0}|}\bigg)^{\frac{p_N-1}{p_N}} \leqslant \gamma \xi \omega_0+2 \gamma \omega_0\bigg(\frac{|Q_{\rho_0}\cap\{u\geqslant \xi \omega_0\}|}{|Q_{\rho_0}|}\bigg)^{\frac{p_N-1}{p_N}}  \ . 
\end{eqnarray*}
 We thereby derive \eqref{eq6.6} once we fix  
 \[\xi=\frac{1}{8\gamma}, \qquad \text{and} \qquad  \alpha=\left( \frac{1}{16\gamma}\right)^{\frac{p_N}{p_N-1}}\,.\]
\end{proof}

\noindent \textbf{Claim $2$.} There exists a time level $s\in (t_0- \rho^{p}_0, t_0-\frac{\alpha}{2} \rho^{p}_0)$ such that
\begin{equation}\label{eq6.7}
\big|\K_{\rho_0}(x_0)\cap\big\{u(\cdot, s) \leqslant \xi \omega_0\big\}\big|\leqslant \frac{1-\alpha}{1-\alpha/2}|\K_{\rho_0}(x_0)|.
\end{equation}

\begin{proof}
The proof is quite simple and is conducted by means of a contradiction argument. Consider that, for all $t\in (t_0- \rho^{p}_0, t_0-\frac{\alpha}{2} \rho^{p}_0)$, inequality \eqref{eq6.7} is violated. Since \eqref{eq6.6} holds true,
\begin{multline*}
(1-\alpha)|Q_{\rho_0}|<\big|\K_{\rho_0}(x_0)\times (t_0- \rho^{p}_0, t_0-\frac{\alpha}{2} \rho^{p}_0)\cap\big\{u\leqslant \xi\,\omega_0\big\}\big|\leqslant\\\leqslant \big|Q_{\rho_0}\cap\big\{u\leqslant \xi \omega_0\big\}\big|\leqslant (1-\alpha)|Q_{\rho_0}|,
\end{multline*}
reaching to a contradiction.
\end{proof}

\vskip0.2cm 

\noindent Now we use Corollary \ref{cor4.1}, for $A=\frac{\alpha}{2}$ and $\bar{\eta}\in (0, 1)$ a constant depending only on the data, to obtain the following inferior bound  
$$u(x, t)\geqslant \bar{\eta}\,\xi \omega_0,\quad (x, t)\in \K_{\rho_0}(x_0)\times(t_0-\frac{\alpha}{4}\rho_0^{p}, t_0),$$
and then, recalling \eqref{eq6.2}, we arrive at the reduction of the oscillation for this second case
\begin{equation}\label{eq6.8}
\osc\limits_{K_{\rho_0}(x_0)\times(t_0-\frac{\alpha}{4}\rho^{p}_0, t_0)} u\leqslant \mu^+_0-\bar{\eta}\xi \omega_0\leqslant(1-\bar{\eta}\xi+\eta)\omega_0=(1-\frac{1}{2}\bar{\eta}\xi)\omega_0,
\end{equation}
provided that 
$$\eta=\frac{1}{2}\bar{\eta}\xi.$$
Combining both reduction of the oscillations \eqref{eq6.5} and \eqref{eq6.8} we arrive at
\begin{equation}\label{eq6.9}
\osc\limits_{Q_{\epsilon \rho_0}}u\leqslant \sigma \omega_0:=\omega_1,\quad \epsilon=\Big(\frac{\alpha}{4}\Big)^{\frac{1}{p}},\quad \sigma=1-\frac{1}{2}\bar{\eta}\xi.
\end{equation}
This proves \eqref{eq6.1} for $j=1$. We can repeat this procedure and prove \eqref{eq6.1}, assuming on each step that
\begin{equation}\label{eq6.10}
\mu^+_j\leqslant (1+\frac{1}{2} \bar{\eta}\xi)\,\omega_j \ , 
\end{equation}
where  $\epsilon$ and $\sigma\in (0, 1)$ are as in \eqref{eq6.9}.  

\vskip0.2cm 

\noindent Assume now that $j_0\geqslant 1$ is the first number such that \eqref{eq6.10} is violated, i.e.
\begin{equation}\label{eq6.11}
\mu^+_{j_0}\geqslant (1+\frac{1}{2} \bar{\eta}\xi)\,\omega_{j_0},\quad \text{and hence}\quad \mu^-_{j_0}\geqslant \bigg(\frac{\bar{\eta}\xi}{2+\bar{\eta}\xi}\bigg) \mu^+_{j_0}.
\end{equation}
In this case, by the fact that \eqref{eq6.10} holds for $j=j_0-1$ we obtain
\begin{equation*}
\bigg(1+\frac{1}{2}\bar{\eta}\xi\bigg)\omega_{j_0}\leqslant \mu^+_{j_0}\leqslant \mu^+_{j_0-1}\leqslant \bigg(1+\frac{1}{2}\bar{\eta}\xi \bigg)\omega_{j_0-1}=\bigg(\frac{2+\bar{\eta}\xi}{2 \sigma}\bigg)\omega_{j_0},
\end{equation*}
and 
\begin{equation}\label{eq6.12}
\frac{1}{2}\bar{\eta}\xi\,\, \omega_{j_0}\leqslant \bigg(\frac{\bar{\eta}\xi}{2+\bar{\eta}\xi}\bigg) \mu^+_{j_0}\leqslant \mu^-_{j_0}\leqslant \mu^+_{j_0}\leqslant \bigg(\frac{1+\bar{\eta}\xi}{2 \sigma} \bigg)\omega_{j_0}.
\end{equation}
Teunction $v:= u/\mu^-_{j_0}$  is a solution of the anisotropic $p$-Laplace evolution equation
\begin{equation}\label{eq6.13}
v_t-\sum_i D_i a_i(x, t, v , D v)=0,\quad (x,t)\in Q_{\rho_{j_0}},
\end{equation}
where $a_i(x, t, v, D v):=\frac{1}{\mu^-_{j_0}}\,A_i(x, t, \mu^-_{j_0} v, \mu^-_{j_0}D v)$
and the following inequalities hold
\begin{equation}\label{eq6.14}
\begin{cases}
a_i(x, t, v, D v) D_i v\geqslant \bar{K}_1(\sigma, \bar{\eta}, \xi) \sum_i |D_i v|^{p_i},\\
|a_i(x, t, v, D v)|\leqslant \bar{K}_2(\sigma, \bar{\eta}, \xi)\, |D_i v|^{p_i-1},\quad i=1, ..., N,
\end{cases}
\end{equation}
with positive constants $\bar{K}_1(\sigma, \bar{\eta}, \xi)$, $\bar{K}_2(\sigma, \bar{\eta}, \xi)$ depending only on $\sigma$, $\bar{\eta}$, $\xi$ and on the constants $K_1$, $K_2$ appearing in \eqref{eq1.2}.

\vskip0.2cm 

\noindent To proceed with the study of the local H\"older continuity of $u$ we need to study separately the degenerate and the singular cases. In both studies we will assume that
\[p_N-p_{1}\leqslant \epsilon_*\, ,\]
for a sufficiently small $\epsilon_*\in (0, 1)$ as in Lemma \ref{lem2.8}.

\subsection{Degenerate Case - $p_N>2$}
Consider the real positive numbers $C_*>1$, $b_2>0$ and $\epsilon_1 \in (0, 1)$, to be defined, and introduce the following numbers
\begin{equation}\label{eq6.15}
\tilde{\mu}^+_{j_0}:=\frac{\mu^+_{j_0}}{\mu^-_{j_0}},\quad \tilde{\mu}^-_{j_0}:=1,\quad \tilde{\omega}_{j_0}:=\frac{\omega_{j_0}}{\mu^-_{j_0}},
\end{equation}
and cylinders
\begin{equation*}
\tilde{Q}^{(\tilde{\omega}_{j_0})}_{\epsilon_1 \rho_{j_0}}:= \K^{\tilde{\omega}_{j_0}/C_*}_{\epsilon_1 \rho_{j_0}}(x_0)\times(t_0-b_2(\epsilon_1 \rho_{j_0})^{p} \tilde{\omega}_{j_0}^{2-p_N}, t_0).
\end{equation*}
We choose $\epsilon_1$ so small that $\tilde{Q}^{(\tilde{\omega}_{j_0})}_{\epsilon_1 \rho_{j_0}}\subset Q_{\rho_{j_0}}$, i.e., by 
\eqref{eq6.12}
\begin{equation*}
\epsilon_1 \Big(C_*\frac{1+\frac{1}{2}\bar{\eta}\xi}{\sigma}\Big)^{\frac{p_N-p_{1}}{p}}\leqslant 1,\quad \text{and}\quad \epsilon^{p}_1\,b_2\,\Big(\frac{1+\frac{1}{2}\bar{\eta}\xi}{\sigma}\Big)^{p_N-2}\leqslant 1.
\end{equation*}
The following  alternative cases are possible: 
either
$$
\big|\K^{\tilde{\omega}_{j_0}/4}_{\epsilon_1 \rho_{j_0}}(x_0)\cap\big\{v\big(\cdot, t_0-b_2(\epsilon_1 \rho_{j_0})^{p} \tilde{\omega}^{2-p_N}_{j_0}\big)\geqslant \tilde{\mu}^+_{j_0}-\frac{1}{4}\tilde{\omega}_{j_0}\big\}\big|\leqslant \frac{1}{2}\big|\K^{\tilde{\omega}_{j_0}/4}_{\epsilon_1 \rho_{j_0}}(x_0)\big|,
$$
or
$$
\big|\K^{\tilde{\omega}_{j_0}/4}_{\epsilon_1 \rho_{j_0}}(x_0)\cap\big\{v\big(\cdot, t_0-b_2(\epsilon_1 \rho_{j_0})^{p} \tilde{\omega}^{2-p_N}_{j_0}\big)\leqslant \tilde{\mu}^-_{j_0}+\frac{1}{4}\tilde{\omega}_{j_0}\big\}\big|\leqslant \frac{1}{2}\big|\K^{\tilde{\omega}_{j_0}/4}_{\epsilon_1 \rho_{j_0}}(x_0)\big|.
$$
By Lemma \ref{lem2.8}, there exist $C_*>1$ and $0<b_1<b_2$, depending only on the data, $\sigma$ and $\bar{\eta}$, such that
\begin{equation}\label{eq6.16}
\text{either} \qquad v(x, t)\leqslant \tilde{\mu}^+_{j_0}-\frac{1}{C_*}\,\tilde{\omega}_{j_0},\qquad 
\text{or} \qquad 
v(x, t)\geqslant \tilde{\mu}^-_{j_0}+\frac{1}{C_*}\,\tilde{\omega}_{j_0},
\end{equation}
 for almost every point $(x, t) \in  \K^{\tilde{\omega}_{j_0}/C_*}_{\epsilon_1 \rho_{j_0}}(x_0)\times (t_0-(b_2-b_1)(\epsilon_1 \rho_{j_0})^{p} \tilde{\omega}^{2-p_N}_{j_0}, t_0)$.
 
 \vskip0.2cm 
 
 \noindent Now we set $\bar{\rho}_{j_0}:=\epsilon_1 \rho_{j_0}$, $\sigma_1:=1-\frac{1}{C_*}$, $\tilde{\omega}_{j_0+1}:=\sigma_1 \tilde{\omega}_{j_0}$ and choose $\epsilon_2 \in (0,1)$ such that
\begin{equation*}
\epsilon^{p}_2\leqslant \,\sigma^{p_N-p_{1}}_1 \quad \text{and}\quad \epsilon^{p}_2\leqslant \big(1-\frac{b_1}{b_2}\big) \sigma^{p_N-2}_1
\end{equation*}
which, together with  \eqref{eq6.16}, guarantee that
\begin{equation}\label{eq6.17}
\osc\limits_{\tilde{Q}^{(\tilde{\omega}_{j_0+1})}_{\bar{\rho}_{j_0+1}}} v\leqslant \tilde{\omega}_{j_0+1},\quad \bar{\rho}_{j_0+1}:=\epsilon_2\,\bar{\rho}_{j_0}; 
\end{equation}
and this is the starting point to an iterative scheme. Indeed, take 
$$\bar{\rho}_{j+j_0}:=\epsilon^j_2 \bar{\rho}_{j_0}\quad \text{and}\quad \tilde{\omega}_{j+j_0}:=\sigma^j_1\,\tilde{\omega}_{j_0},\quad j\geqslant 0.$$
Then, starting from \eqref{eq6.17} and applying several times Lemma \ref{lem2.8}, we obtain 
\begin{equation}\label{eq6.18}
\osc\limits_{\tilde{Q}^{(\tilde{\omega}_{j+j_0})}_{\bar{\rho}_{j+j_0}}} v\leqslant \tilde{\omega}_{j+j_0},\quad \text{for all}\quad j\geqslant 0.
\end{equation}
Observe that under condition \eqref{eq6.12} and the  additional assumption
\begin{equation*}
\epsilon^{p}_2\leqslant \Big(\frac{\bar{\eta} \xi}{2}\Big)^{p_N-p_{1}}\quad \text{and}\quad \epsilon^{p}_2\leqslant \Big(\frac{\bar{\eta} \xi}{2}\Big)^{p_N-2},
\end{equation*}

$$Q_{\bar{\rho}_{j+1+j_0}}\subset \tilde{Q}^{(\tilde{\omega}_{j+j_0})}_{\bar{\rho}_{j+j_0}},\quad \text{for all}\quad j\geqslant 0, $$
and hence, redefining $\bar{\rho}_{j+j_0}$ and $\omega_{j+j_0}$ as $\bar{\rho}^{(1)}_{j+j_0}:=\epsilon_2 \bar{\rho}_{j+j_0}=\epsilon^{j+1}_2 \bar{\rho}_{j_0}$, $\omega_{j+j_0}:=\sigma^j_1 \omega_{j_0}$, inequality \eqref{eq6.18} yields 
\begin{equation}\label{eq6.19}
\osc\limits_{Q_{\bar{\rho}^{(1)}_{j+j_0}}} u=\osc\limits_{Q_{\bar{\rho}_{j+1+j_0}}} u\leqslant \omega_{j+j_0},\quad \text{for all}\quad j\geqslant 0,
\end{equation}
which completes the construction of the sequences as in \eqref{eq6.1} in  the case if \eqref{eq6.11} holds.
Now the H\"{o}lder continuity follows by standard iterative real analysis methods, see, for example \cite[Chap.3]{DiB}.

\subsection{Singular Case - $p_N\leqslant 2$}
Define numbers $\tilde{\mu}^+_{j_0}$, $\tilde{\mu}^-_{j_0}$ and $\tilde{\omega}_{j_0}$ as in \eqref{eq6.15},  and
for fixed numbers $C_*>1$, $b_2>0$, $\epsilon_1 \in (0, 1)$  introduce the cylinders
\begin{equation*}
\tilde{Q}^{(\tilde{\omega}_{j_0})}_{\epsilon_1 \rho_{j_0}}:=
 \tilde{\K}^{(C_*)}_{\epsilon_1 \rho_{j_0}}\times (t_0-b_2(\epsilon_1 \rho_{j_0})^{p}, t_0),
\end{equation*}
where 
\begin{equation*}
\tilde{\K}^{(C_*)}_{\epsilon_1 \rho_{j_0}}:=\{x: |x_i-x_{i,0}|\leqslant C^{\frac{p_N-p_i}{p_i}}_*\tilde{\omega}^{-\frac{2-p_i}{p_i}}_{j_0}(\epsilon_1 \rho_{j_0})^{\frac{p}{p_i}},\quad i=1, ..., N\}.
\end{equation*}
Note that,  redefining $\rho_{j_0}$ by $r^{p}_{j_0}=\tilde{\omega}^{p_N-2}_{j_0}\, \rho^{p}_{j_0}$, then $\tilde{\K}^{(C_*)}_{\epsilon_1 \rho_{j_0}}$ coincides with $\K^{\tilde{\omega}_{j_0}/C_*}_{\epsilon_1 r_{j_0}}$ 
and the cylinder $\tilde{Q}^{(\tilde{\omega}_{j_0})}_{\epsilon_1 \rho_{j_0}}$ coincides with the cylinder $\tilde{Q}^{(\tilde{\omega}_{j_0})}_{\epsilon_1 r_{j_0}}$ defined in the previous section.

\vskip0.2cm 

\noindent Choose $\epsilon_1$ small enough such that $\tilde{Q}^{(\tilde{\omega}_{j_0})}_{\epsilon_1 \rho_{j_0}}\subset Q_{\rho_{j_0}}$, i.e., by 
\eqref{eq6.12}
\begin{equation*}
\epsilon_1 C^{\frac{p_N-p_{1}}{p}}_* \Big(\frac{2+\bar{\eta}\xi}{2 \sigma}\Big)^{\frac{2-p_{1}}{p}}\leqslant 1 \quad \text{and}\quad \epsilon^{p}_1 b_2 \leqslant 1.
\end{equation*}
We distinguish  two different cases: either
\begin{equation*}
\big|\tilde{\K}^{(4)}_{\epsilon_1 \rho_{j_0}}\cap\big\{v\big(\cdot, t_0-b_2(\epsilon_1 \rho_{j_0})^{p}\big)\geqslant \tilde{\mu}^+_{j_0}-\frac{1}{4}\omega_{j_0}\big\}\big|\leqslant \frac{1}{2}|\tilde{\K}^{(4)}_{\epsilon_1 \rho_{j_0}}|,
\end{equation*}
or
\begin{equation*}
\big|\tilde{\K}^{(4)}_{\epsilon_1 \rho_{j_0}}\cap\big\{v\big(\cdot, t_0-b_2(\epsilon_1 \rho_{j_0})^{p}\big)\leqslant \tilde{\mu}^-_{j_0}+\frac{1}{4}\tilde{\omega}_{j_0}\big\}\big|\leqslant \frac{1}{2}|\tilde{\K}^{(4)}_{\epsilon_1 \rho_{j_0}}|.
\end{equation*}
 Lemma \ref{lem2.8}  yields the existence of $C_*>1$ and $0<b_1<b_2$, depending only on the data, $\sigma$ and $\bar{\eta}$, such that
\begin{equation}\label{eq6.20}
\text{either} \qquad v(x, t)\leqslant \tilde{\mu}^+_{j_0}-\frac{1}{C_*}\,\tilde{\omega}_{j_0},\qquad 
\text{or}\qquad v(x, t)\geqslant \tilde{\mu}^-_{j_0}+\frac{1}{C_*}\,\tilde{\omega}_{j_0},
\end{equation}
for almost every $(x, t) \in \tilde{\K}^{(C_*)}_{ \epsilon_1 \rho_{j_0}}\times (t_0-(b_2-b_1)(\epsilon_1 \rho_{j_0})^{p}, t_0)$.

\vskip0.2cm 

\noindent Set $\bar{\rho}_{j_0}:=\epsilon_1 \rho_{j_0}$, $\sigma_1:=1-\frac{1}{C_*}$, $\tilde{\omega}_{j_0+1}:=\sigma_1 \tilde{\omega}_{j_0}$ and choose $\epsilon_2 \in (0,1)$ such that
\begin{equation*}
\epsilon^{p}_2\leqslant\,\sigma^{2-p_{1}}_1,\quad \text{and}\quad \epsilon^{p}_2\leqslant \big(1-\frac{b_1}{b_2}\big).
\end{equation*}
Then, by \eqref{eq6.20}, one obtains
\begin{equation}\label{eq6.21}
\osc\limits_{\tilde{Q}^{(\tilde{\omega}_{j_0+1})}_{\bar{\rho}_{j_0+1}}} v\leqslant \tilde{\omega}_{j_0+1},\quad \bar{\rho}_{j_0+1}:=\epsilon_2\,\bar{\rho}_{j_0}.
\end{equation}
We're just one step away from presenting the final part of the iterative scheme leading to H\"older's continuity. Let's first choose 
$$\bar{\rho}_{j+j_0}:=\epsilon^j_2 \bar{\rho}_{j_0}\quad \text{and}\quad \tilde{\omega}_{j+j_0}:=\sigma^j_1\,\tilde{\omega}_{j_0},\quad j\geqslant 0,$$
then we repeat the previous procedure and 
\begin{equation}\label{eq6.22}
\osc\limits_{\tilde{Q}^{(\tilde{\omega}_{j+j_0})}_{\bar{\rho}_{j+j_0}}} v\leqslant \tilde{\omega}_{j+j_0},\quad \text{for all}\quad j\geqslant 0,
\end{equation}
and finally, by considering the additional assumption on $\epsilon_2$ 
\begin{equation*}
\epsilon^{p}_2\leqslant \big(\frac{\bar{\eta}\xi}{2}\big)^{2-p_{1}} \quad \text{and}\quad \epsilon^{p}_2\leqslant b_2,
\end{equation*}
we obtain that
$$Q_{\bar{\rho}_{j+1+j_0}}\subset \tilde{Q}^{(\tilde{\omega}_{j+j_0})}_{\bar{\rho}_{j+j_0}},\quad \text{for all}\quad j\geqslant 0. $$
Redefine  $\bar{\rho}_{j+j_0}$ and $\omega_{j+j_)}$ as $\bar{\rho}^{(1)}_{j+j_0}:=\epsilon_2 \bar{\rho}_{j+j_0}=\epsilon^{j+1}_2 \bar{\rho}_{j_0}$ and $\omega_{j+j_0}:=\sigma^j_1 \omega_{j_0}$, then  \eqref{eq6.22} yields 
\begin{equation}\label{eq6.23}
\osc\limits_{Q_{\bar{\rho}^{(1)}_{j+j_0}}} u=\osc\limits_{Q_{\bar{\rho}_{j+1+j_0}}} u\leqslant \omega_{j+j_0},\quad \text{for any}\quad j\geqslant 0,
\end{equation}
which completes the construction of the sequences as in \eqref{eq6.1} in  the case if \eqref{eq6.11} holds. Once again and as in the degenerate case, the H\"{o}lder continuity follows by standard iterative real analysis methods. This completes the proof of Theorem \ref{th1.2}.

\section{Appendix}\label{appendix}

\subsection{Motivation for the Functional Setting} \label{redefinition}\noindent The functional assumption \eqref{functional} is natural for the integrals of \eqref{eq1.3} to be convergent. Indeed, by by H\"older inequality, for any $0<t<T$, we can estimate the first parabolic term as
\[\int\limits_E u(x,t)\,\varphi(x,t)\,dx \leq \gamma \bigg(\sup_{[0,T]} \|u(\cdot, t)\|_{L^{\frac{p_N}{p_N-1}}(E)} \bigg) \bigg(\sup_{[0,T]}\|\varphi(\cdot, t)\|_{L^{p_N}(E)}\bigg)< \infty\] since $u \in C_{loc}(0,T; L^{\frac{p_N}{p_N-1}}(E)$ and $\varphi \in C_{loc}(0,T; L^{p_N}(E)) \subset W^{1,p_N}_{loc}(0,T; L^{p_N}(E))$. Similarly we bound the second parabolic term, thanks to the previous assumption since 
\begin{equation} \label{si} u \in C_{loc}(0,T; L^{\frac{p_N}{p_N-1}}_{loc}(\Omega))\subset L^{\frac{p_N}{p_N-1}}_{loc}(\Omega_T),\end{equation} with
\[\int_{t_1}^{t_2} \int_E u \partial_t \varphi \, dxdt \leq \gamma \bigg( \int_{t_1}^{t_2} \int_E u^{\frac{p_N}{p_N-1}} \bigg)^{\frac{p_N-1}{p_N}} \bigg( \int_{t_1}^{t_2} \int_E |\partial_t \varphi|^{p_N} \bigg)^{\frac{1}{p_N}}< \infty\,.\]
Finally, the elliptic term is controlled by means of the structure conditions \eqref{eq1.2}. For all $i \in \{1, \dots,N\}$
\begin{equation*}
\begin{aligned}
    \int\limits^{t_2}_{t_1}\int\limits_E A_i(x, t, u, D u)\,D_i \varphi\, \, dx dt  \leq \gamma \| D_i(u^{\frac{1}{p_i-1}}) \|_{L^{p_i-1}(E\times (t_1,t_2)} \, \|D_i \varphi\|^{p_i}_{L^{p_i}(E\times (t_1,t_2))}\,,
\end{aligned}
\end{equation*} showing therefore the necessity of  the assumption $D_i u^{\frac{1}{p_i-1}} \in L^{p_i}_{loc}(0,T; L^{p_i}_{loc}(E))$. Finally, the  assumption \eqref{si} allows similarly for the finiteness of the energy inequalities \eqref{eq2.3}, and will be actually used crucially to prove them here below in Section \ref{proof-energy-estimates}.

\subsubsection{Existence Theory} Boundary value problems are well-posed (see \cite{Vestberg}) when considering equation \eqref{eq1.1}-\eqref{eq1.2} within the special form of 
\[A_i(x,t,u,Du)=a_{i,j}(x,t,u) |D_i(u^{\frac{1}{p_i-1}})|^{p_i-2}D_i(u^{\frac{1}{p_i-1}})\,,\] a boundary datum $g \in L^{\bf{p}}(0,T; W^{1,\bf{p}}(\Omega)) \cap L^{\infty}(\Omega_T)$ such that $\partial_t g \in L^2(\Omega_T)$, together with an initial value $0 \leq u_0 \in L^{\infty}(\Omega)$. We refer to the study \cite{Vestberg} with $m_i:= 1/(p_i-1)$. Two are the conditions required on the spread of $p_i$s: the first one is
\[1<p_i<2 \qquad \forall i=1,\dots, N\,,\]
which is nevertheless claimed to be unnecessary in \cite{Vestberg}; and the condition 
\[\frac{1}{p_i-1}\leq \frac{p_i}{p_i-1}\bigg( \frac{1}{p_N-1}\bigg) \quad \iff \quad p_N<p_1+1\,.\]
This last condition is implied by the assumption $p_N-p_1<\epsilon$ of Theorem \ref{th1.2} and in particular it provides, for the solutions, the existence of a spatial gradient of $D(u^{\frac{1}{p_N-1}})$ together with the chain rule, see Section 4 of \cite{Vestberg}.

\subsubsection{The structure conditions revisited}
In order to end up with positive exponents in the energy estimates, we rewrite the structure conditions \eqref{eq1.2} as follows

\begin{equation}\label{no}
\begin{cases}
 A_i(x, t, u, D u)\,D_i u^{\frac{1}{p_N-1}}  \geqslant  \tilde{K}_{1}\, u^{\frac{p_N-p_i}{p_N-1}}\,|D_i u^{\frac{1}{p_N-1}}|^{p_i},\\ $\quad$ \\
|A_i(x, t, u,  D u)| \leqslant  \tilde{K}_{2}\,  u^{\frac{p_N-p_i}{p_N-1}}|D_i u^{\frac{1}{p_N-1}}|^{p_i-1},\qquad  \forall \, i=1, ...,N,
\end{cases}
\end{equation}
where $\tilde{K}_1$, $\tilde{K}_2$ are positive constants,  and $Du = (D_1 u,... , D_N u)$ is assumed (by the request \eqref{functional}) to satisfy weakly the chain rule  
\begin{equation}
    \label{chain-rule}D_i u = (p_N-1)u^{\frac{p_N-2}{p_N-1}}\,D_i \big(u^{\frac{1}{p_N-1}}\big), \quad\qquad \quad \forall i=1, ...,N.
\end{equation} The fact that \eqref{no} implies \eqref{eq1.2} under the assumption of the chain rule, can be seen by simple algebraic computations as
\[ \bigg(\frac{u^{{2-p_N}}}{\gamma} \bigg) A_i \,  D_i u  =A_i\,  D_i (u^{\frac{1}{p_N-1}})  \ge  \tilde{K}_1 u^{\frac{p_N-p_i}{p_N-1}}|D_i(u^{\frac{1}{p_N-1}})|^{p_i} = \frac{\tilde{K}_1}{\gamma} u^{\frac{p_N-p_i}{p_N-1}+(\frac{2-p_N}{p_N-1})p_i} |D_i u |^{p_i}\,.  \]
Manipulating the terms in the set $[u>0]\cap \Omega_T$, we obtain the first condition in \eqref{eq1.2}; and analogously for the second one.

\subsection{Alternative weak formulation} \label{mollificazione}
\noindent Since we deal with measurable and bounded coefficients, the time derivative of a weak solution to \eqref{eq1.3} exists only in the sense of distribution (see for an example the diffraction problem in \cite{LadSolUra}, Sec. 13 of Ch. III, pages 224-233). In 
order to overcome the lack of regularity in the time variable, we define the following 
mollification in time: let $v \in L^1(E \times (t_1,t_2))$ for a compact set $E \subset \Omega$, times $0<t_1<t_2<T$, and let us set 
\[\llbracket v \rrbracket_h(x, t):=\frac{1}{h}\int\limits_{t_1}^{t} e^{\frac{\tau-t}{h}}\,v(x, \tau)\,d\tau,\quad \qquad \llbracket v \rrbracket_{\bar{h}}(x, t):=\frac{1}{h}\int\limits_t^{t_2} e^{\frac{t-\tau}{h}}\,v(x, \tau)\,d\tau,\]
We refer the reader to \cite{KinLin} and the Appendix B of \cite{BDM} for the proofs of the following properties.
\begin{lemma}\label{lem-moll}
For any $s\geqslant 1$ we have
\vskip0.2cm \noindent 
$(i)$ If $v\in L^s(\Omega_T)$, then $\llbracket v \rrbracket_h\in L^s(\Omega_T)$. Moreover, $||\llbracket v \rrbracket_h||_{L^s(\Omega_T)}
\leqslant ||v||_{L^s(\Omega_T)}$, and $\llbracket v \rrbracket_h \rightarrow v$ in $L^s(\Omega_T)$ and almost everywhere on $\Omega_T$
as $h\rightarrow 0$.
\vskip0.2cm \noindent 
$(ii)$ $\llbracket v \rrbracket_h \in C(0, T; L^s(\Omega))$.
\vskip0.2cm \noindent 
$(iii)$ Almost everywhere on $\Omega_T$ there hold
$$\frac{\partial}{\partial t} \llbracket v \rrbracket_h =\frac{1}{h}\big(v-\llbracket v \rrbracket_h\big).$$
\vskip0.2cm \noindent 
$(iV)$ If $Dv \in L^s(\Omega_T)$, then $D \llbracket v \rrbracket_h = \llbracket D v \rrbracket_h \rightarrow D v$ in $L^s(\Omega_T)$
and almost everywhere on $\Omega_T$ as $h \rightarrow 0$.
\vskip0.2cm \noindent 
$(V)$ If $v\in C(0, T; L^s(\Omega))$, then $\llbracket v \rrbracket_h(\cdot, t) \rightarrow v(\cdot, t)$ in $L^s(\Omega)$ and almost everywhere on $\Omega$ for every $t\in (0, T)$ as $h\rightarrow 0$.
\end{lemma}

\noindent From the  differential 
inequalities \eqref{eq1.3}, we deduce the mollified version as in Lemma 4.2 of \cite{Vestberg}, 
\begin{equation}\label{eq.mollified}
\int_{t_1}^{t_2} \int_E \Big(\frac{\partial}{\partial t} \llbracket u \rrbracket_h\,\varphi+\sum_i \llbracket A_i(x, \tau, u, D u) \rrbracket_h D_i \varphi\Big)\, dx d\tau \leqslant (\geqslant)\, \int_E u(x,t_1) \llbracket\varphi\rrbracket_{\bar{h}} (x,t_1) \, dx\,,
\end{equation}
for any non-negative $\varphi\in L^{\textbf{p}}(0, T; W^{1, \textbf{p}}_0(\Omega))$.

\subsection{Energy Estimates. Proof of Lemma 3.6}
\label{proof-energy-estimates} We only consider the case of a local weak sub-solution to \eqref{eq1.3}, the case of super-solutions can be considered completely analogously.\vskip0.1cm \noindent For
fixed $\tau-\eta<t_1<t_2<\tau$ and sufficiently small $\epsilon \in (0,1)$ we define (see Figure \ref{fig:placeholder})
\begin{equation*}
\psi_\epsilon(t):=
\begin{cases}
0,\quad \tau-\eta\leqslant t\leqslant t_1-\epsilon,\\
1+\frac{t-t_1}{\epsilon},\quad t_1-\epsilon\leqslant t \leqslant t_1,\\
1,\quad t_1\leqslant t \leqslant t_2,\\
1-\frac{t-t_2}{\epsilon},\quad t_2\leqslant t\leqslant t_2+\epsilon,\\
0,\quad t_2+\epsilon \leqslant t\leqslant \tau,
\end{cases}
\end{equation*}
and choose in \eqref{eq.mollified} the test function $\varphi= (u^{\frac{1}{p_N-1}}-k^{\frac{1}{p_N-1}})_+\,\zeta^{p_N}(x, t)\,\psi_\epsilon(t).$ So we obtain

\begin{figure}
    \centering
  
\begin{tikzpicture}[x=1cm,y=1cm,>=stealth]
  \def\tau{9}          
  \def\eta{7}          
  \def\tone{3}         
  \def\ttwo{6}         
  \def\eps{0.8}        
  \pgfmathsetmacro{\taumenoeta}{\tau-\eta}
  \pgfmathsetmacro{\tunoLess}{\tone-\eps}
  \pgfmathsetmacro{\tduePlus}{\ttwo+\eps}
  \draw[->] (0,0) -- (\tau+0.8,0) node[below] {$t$};
  \draw[->] (0,-0.1) -- (0,1.4) node[left] {$\textcolor{blue}{\psi_\epsilon(t)}$};
  \draw[thick,blue] 
        (0,0) --
        (\taumenoeta,0) --
        (\tunoLess,0)   --
        (\tone,1)       --
        (\ttwo,1)       --
        (\tduePlus,0)   --
        (\tau,0);
  \foreach \x/\label in {
                         \tunoLess/{$t_1-\epsilon$},%
                         \tone/{$t_1$},%
                         \ttwo/{$t_2$},%
                         \tduePlus/{$t_2+\epsilon$},%
                         }
     {
       \draw[dashed] (\x,0) -- (\x,-0.1);
       \node[below] at (\x,-0.1) {\label};
     }
  \draw[dashed] (0,1) -- (-0.1,1);
  \node[left] at (-0.1,1) {$1$};
  \foreach \X/\Y in {(\tunoLess,0),(\tone,1),(\ttwo,1),(\tduePlus,0)}
      \fill \X \Y circle (1.2pt);
\end{tikzpicture}

    \caption{Illustration of the trapezoidal function $\psi_{\epsilon}$.}
    \label{fig:placeholder}
\end{figure}

\begin{equation}\label{ciccio}
\begin{aligned}
     0 \ge \, I_1+I_2-I_3&= \int_{\tau-\eta}^{\tau} \int_E  \partial_t \llbracket u \rrbracket_{h} (u^{\frac{1}{p_N-1}}-k^{\frac{1}{p_N-1}})_+ \zeta^{p_N} \psi_{\epsilon} dxdt +\\
    \\ & \qquad \quad + \sum_i \int_{\tau-\eta}^{\tau} \int_E \llbracket A_i \rrbracket_h  \, D_i \bigg[ (u^{\frac{1}{p_N-1}}-k^{\frac{1}{p_N-1}})_+\,\zeta^{p_N}\,\psi_\epsilon\bigg] dxdt-\\
    & \qquad \qquad  \qquad\qquad  \qquad -\int_E u(x,\tau-\eta) \bigg(\int_{\tau-\eta}^{\tau} \frac{e^{\frac{\tau-\eta-s}{h}}}{h} \varphi(x, s) ds \bigg) dx \,.
\end{aligned}
\end{equation}
\noindent Our intention is to pass to the limit $h \downarrow 0$. The last integral vanishes, thanks to the vanishing property $\varphi(x,t)=0$ for $t \in [\tau-\eta\, , t_1-\epsilon]$, the time-continuity {\it à la Bochner}, Fubini-Tonelli's theorem and elementary integration, since

\begin{equation*}
    \begin{aligned}
\int_E u(x,\tau-\eta)& \bigg(\int_{\tau-\eta}^{\tau} \frac{e^{\frac{\tau-\eta-s}{h}}}{h} \varphi(x, s) ds \bigg) dx = \frac{1}{h}\int_{t_1-\epsilon}^{\tau} e^{\frac{\tau-\eta-s}{h}} \bigg( \int_E u(x,\tau-\eta)\varphi(x,s)\, dx \bigg) d s\\
& \leq  \bigg(\sup_{t_1-\epsilon<\tau<t_2} \|u(x,t_1)\varphi(x,\tau)\|_{L^1(E)} \bigg) (e^{\frac{\tau-\eta - t_1+\epsilon}{h}}) \rightarrow 0\,,
    \end{aligned}
\end{equation*} when $h$ vanishes. For $I_1$ we compute
\begin{multline}\label{1.1}
\int\limits_{\tau-\eta}^{\tau}\int\limits_{E}\,\frac{\partial}{\partial t} \llbracket u \rrbracket_h\,\varphi\,dx dt=
\int\limits_{\tau-\eta}^{\tau}\int\limits_{E}\,\frac{\partial}{\partial t} \llbracket u \rrbracket_h\,\big(\llbracket u\rrbracket^{\frac{1}{p_N-1}}_h -k^{\frac{1}{p_N-1}}\big)_+\,
\zeta^{p_N}\,\psi_\epsilon dx dt+\\+\int\limits_{\tau-\eta}^{\tau}\int\limits_{E}\,\frac{\partial}{\partial t} \llbracket u \rrbracket_h\,\Big[(u^{\frac{1}{p_N-1}}-k^{\frac{1}{p_N-1}})_+-\big(\llbracket u\rrbracket^{\frac{1}{p_N-1}}_h -k^{\frac{1}{p_N-1}}\big)_+\Big]\,
\zeta^{p_N}\,\psi_\epsilon dx dt= I_{1,1}+I_{1,2}.
\end{multline}
To estimate $I_{1,1}$ we use the properties of $\psi_{\epsilon}$  to get
\begin{equation*}
\begin{aligned}
I_{1,1}=&\int\limits_{\tau-\eta}^{\tau}\int\limits_{E}\,\frac{\partial}{\partial t}\bigg[\int\limits^{\llbracket u\rrbracket_h}_{k}(z^{\frac{1}{p_N-1}}-k^{\frac{1}{p_N-1}})_+\,dz \,\,\zeta^{p_N}\,\psi_\epsilon\bigg] dxdt-\\
&- p_N\int\limits_{\tau-\eta}^{\tau}\int\limits_{E}\int\limits^{\llbracket u\rrbracket_h}_{k}(z^{\frac{1}{p_N-1}}-k^{\frac{1}{p_N-1}})_+\,dz 
\,\,\zeta^{p_N-1}\zeta_t\,\psi_\epsilon dx dt-\\
&\qquad -\int\limits_{\tau-\eta}^{\tau}\int\limits_{E}\int\limits^{\llbracket u\rrbracket_h}_{k}(z^{\frac{1}{p_N-1}}-k^{\frac{1}{p_N-1}})_+\,dz 
\,\,\zeta^{p_N}\,\psi'_\epsilon dx dt\geqslant\\
&\qquad \qquad \geqslant -p_N\int\limits^{t_2+\epsilon}_{t_1-\epsilon}\int\limits_{E}\int\limits^{\llbracket u\rrbracket_h}_{k}(z^{\frac{1}{p_N-1}}-k^{\frac{1}{p_N-1}})_+\,dz 
\,\,\zeta^{p_N-1}|\zeta_t|\,dx dt-\\
&\qquad \qquad \qquad \qquad -\int\limits^{t_2+\epsilon}_{t_1-\epsilon}\int\limits_{E}\int\limits^{\llbracket u\rrbracket_h}_{k}(z^{\frac{1}{p_N-1}}-k^{\frac{1}{p_N-1}})_+\,dz 
\,\,\zeta^{p_N}\,\psi'_\epsilon dx dt.
\end{aligned}
\end{equation*}
On the other hand, using property $(iii)$ from Lemma \ref{lem-moll} and the fact that the function $s \rightarrow \max\{s^\frac{1}{p_N-1}, k^\frac{1}{p_N-1}\}$ is monotone increasing, we see that $I_{1,2}$ can be discarded since
\begin{equation*}
I_{1,2}=\frac{1}{h}\int\limits_{\tau-\eta}^{\tau}\int\limits_{E}\,\Big[u-\llbracket u\rrbracket_h\Big]\Big[(u^{\frac{1}{p_N-1}}-k^{\frac{1}{p_N-1}})_+-\big(\llbracket u\rrbracket^{\frac{1}{p_N-1}}_h -k^{\frac{1}{p_N-1}}\big)_+\Big]\,
\zeta^{p_N}\,\psi_\epsilon dx dt\geqslant 0.
\end{equation*}
We can let $h\rightarrow 0$ on $I_1$ and invoke Lebesgue's theorem of dominated convergence relying on the assumption \[u^{\frac{1}{p_N-1}} \in L^1_{loc}(\Omega_T) \subset L^{\bf{p}}_{loc}(0,T;W^{1,\bf{p}}_{loc}(\Omega))\] and property $(i)$ of Lemma \ref{lem-moll}. By using the fact that $g_{+}(u^{\frac{1}{p_N-1}}, k^{\frac{1}{p_N-1}})=\int\limits^{u}_{k}(z^{\frac{1}{p_N-1}}-k^{\frac{1}{p_N-1}})_+ dz$ is continuous w.r.t $u$, the limit can be computed
\begin{multline*}
\lim\limits_{h\rightarrow 0}\int\limits_{\tau-\eta}^{\tau}\int\limits_{E}\,\frac{\partial}{\partial t} \llbracket u \rrbracket_h\,\varphi\,dx dt\geqslant
-p_N\int\limits^{t_2+\epsilon}_{t_1-\epsilon}\int\limits_{E}g_+(u^{\frac{1}{p_N-1}}, k^{\frac{1}{p_N-1}}) 
\,\,\zeta^{p_N-1}|\zeta_t|\, dx dt-\\-\int\limits^{t_2+\epsilon}_{t_1-\epsilon}\int\limits_{E} g_+(u^{\frac{1}{p_N-1}}, k^{\frac{1}{p_N-1}})
\,\,\zeta^{p_N}\,\psi'_\epsilon dx dt.
\end{multline*}
Now we can pass to the limit $\epsilon \rightarrow 0$ in this inequality and obtain
\begin{multline*}
\lim\limits_{\epsilon \rightarrow 0}\lim\limits_{h\rightarrow 0}\int\limits_{\tau-\eta}^{\tau}\int\limits_{E}\,\frac{\partial}{\partial t} \llbracket u \rrbracket_h\,\varphi\,dx dt\geqslant\int\limits_{E} g_+(u^{\frac{1}{p_N-1}}, k^{\frac{1}{p_N-1}})
\,\,\zeta^{p_N}\, dx\bigg|^{t_2}_{t_1}-\\-p_N\int\limits^{t_2}_{t_1}\int\limits_{E}g_+(u^{\frac{1}{p_N-1}}, k^{\frac{1}{p_N-1}}) 
\,\,\zeta^{p_N-1}|\zeta_t|\, dx dt.
\end{multline*}
In order to pass to the limit for $h\downarrow 0$ the diffusion term $I_2$, we split it by the product of weak derivatives as
\begin{equation*}
    \begin{aligned}
        I_2= \int\limits_{\tau-\eta}^{\tau}\int\limits_{E}& \llbracket A_i \rrbracket_h  \, \bigg(D_i (u^{\frac{1}{p_N-1}}-k^{\frac{1}{p_N-1}})_+\bigg)\,\zeta^{p_N}\,\psi_\epsilon dxdt +\\
        &+p_N\int\limits_{\tau-\eta}^{\tau}\int\limits_{E} \llbracket A_i \rrbracket_h  \, (u^{\frac{1}{p_N-1}}-k^{\frac{1}{p_N-1}})_+\,\zeta^{p_N-1} D_i \zeta\,\psi_\epsilon dxdt \\
        &\qquad \qquad  = I_{2,1}+ I_{2,2},
    \end{aligned}
\end{equation*} and apply Lebesgue's dominated convergence theorem to pass to the limit, since by point $(i)$ of Lemma \ref{lem-moll}, structure conditions \eqref{eq1.3} and \eqref{functional} we have
\[ | A_i(x,t,u,Du) | \leq \gamma |D_i u^{\frac{1}{p_i-1}}|^{p_i-1} \quad \Rightarrow \quad \llbracket A_i \rrbracket_h \in L^{\frac{p_i}{p_i-1}}_{loc}(\Omega_T).\]
Hence, the approximated inequality \eqref{ciccio} as soon as both $h$ and $\epsilon$ vanish tends to 
\begin{equation}\label{ciccio2}
\begin{aligned}
     0 \ge& \, \int\limits_{E} g_+(u^{\frac{1}{p_N-1}}, k^{\frac{1}{p_N-1}})
\,\,\zeta^{p_N}\, dx\bigg|^{t_2}_{t_1}-p_N\int\limits^{t_2}_{t_1}\int\limits_{E}g_+(u^{\frac{1}{p_N-1}}, k^{\frac{1}{p_N-1}}) 
\,\,\zeta^{p_N-1}|\zeta_t|\, dx dt\, +\\
&\qquad +\sum_i \int_{\tau-\eta}^{\tau} \int_E  A_i(x,t,u,Du)  \, D_i (u^{\frac{1}{p_N-1}}-k^{\frac{1}{p_N-1}})_+\,\zeta^{p_N}\, dxdt\, +\\
&\qquad \qquad -p_N\sum_i \int_{\tau-\eta}^{\tau} \int_E  |A_i(x,t,u,Du)|  \, (u^{\frac{1}{p_N-1}}-k^{\frac{1}{p_N-1}})_+\,\zeta^{p_N-1} |D_i \zeta|\, dxdt\,.
\end{aligned}
\end{equation}
\noindent Now, the standard use of the structure conditions \eqref{eq1.3} and the Young inequality yield the desired estimate \eqref{eq2.3}.

\subsection{Expansion of Positivity for Anisotropic $p$-Laplace Evolution Equation}
For the proof of H\"{o}lder's continuity of solutions to \eqref{eq1.1}-\eqref{eq1.2} we will make use of the following result which can be found in \cite{CiaHenSavSkr2}. For that purpose, consider the equation
\begin{equation}\label{eq2.5}
u_t-\sum_i D_i a_i(x, t, u, Du)=0, \quad (x,t)\in \Omega_T,
\end{equation}
and let the following inequalities hold
\begin{equation}\label{eq2.6}
\begin{cases}
\sum\limits_{i=1}^N a_i(x, t, u, Du) D_iu\geqslant \bar{K}_1 \sum_i |D_i u|^{p_i},\\
|a_i(x, t, u, D u)|\leqslant \bar{K}_2 \Big(\sum_j |D_j u|^{p_j}\Big)^{1-\frac{1}{p_i}},\quad i=1, ..., N,
\end{cases}
\end{equation}
for fixed positive constants $\bar{K}_1$ and $\bar{K}_2.$  
\vskip0.2cm \noindent For a point $(y,s) \in \Omega_T$, and constants $C_*>1$ and $b$, $r>0$ let $\mu^+, \mu^-, \omega$ be positive numbers satisfying
\[\mu^+\geqslant \sup\limits_{Q^{\omega/C_*}_{r, \theta}(y,s+\theta)}u,\qquad \quad \mu^-\leqslant \inf\limits_{Q^{\omega/C_*}_{r, \theta}(y,s+\theta)}u,\quad  \qquad \omega\geqslant \mu^+-\mu^-\, ,\] and let $\theta$ be an intrinsic time parameter defined as
\[ \theta:=b\,r^p\,\omega^{2-p_N} \ . \]
Finally, let us set
\[v_+:=\mu^+-u,\quad v_-:=u-\mu^-.\]
\noindent 
The expansion of positivity is a property of locally bounded local weak solutions to \eqref{eq2.5}-\eqref{eq2.6}, as shown in Theorems 3.1 and 3.4 of \cite{CiaHenSavSkr2}, at the price of restricting the range $p_N-p_1$ of a qualitative parameter $\epsilon_*$ that is determined only from the data.

\begin{lemma}\label{lem2.8}
Let $u\in C_{loc}(0, T; L^2_{loc}(\Omega))\cap L^{\textbf{p}}_{loc}(0, T; W^{1, \textbf{p}}_{loc}(\Omega))\cap L^\infty_{loc}(\Omega_T)$ be a local weak solution to \eqref{eq2.5}-\eqref{eq2.6} in $\Omega_T$,  and assume also  that
\begin{equation}\label{eq2.7}
\big|\K^{\xi \omega}_r(y)\cap\big\{v_{\pm}(\cdot, s)\geqslant \xi \omega\big\}\big|\geqslant \alpha_0 |\K^{\xi \omega}_r(y)|,
\end{equation}
with some $\xi$, $\alpha_0 \in (0, 1)$. Then there exist  numbers $C_*>1$, $0<\bar{b}<b$, $\epsilon_* \in (0, 1)$ depending only on the data, $\xi$ and $\alpha_0$
such that
\begin{equation}\label{eq2.8}
v_{\pm}(x, t)\geqslant \frac{\omega}{C_*},\quad (x, t)\in \K^{\omega/C_*}_{r}(y)\times (s+\bar{b}r^p \omega^{2-p_N}, s+b r^p \omega^{2-p_N}),
\end{equation}
provided that  $Q^{\omega/C_*}_{8 r, 8\theta}(y,s+\theta)\subset \Omega_T$ and
\begin{equation}\label{eq2.9}
p_N-p_{1}\leqslant \epsilon_*.
\end{equation}
\end{lemma}

\subsection{Auxiliary Lemmas}
The following lemma is the well-known De Giorgi-Poincar\'e lemma (see
\cite[Lemma 2.2, Chap.2]{DiBGiaVes3}), presented for  standard cubes.
\begin{lemma}\label{lem2.3}
Let $u\in W^{1,1}(K_{1}(y))$ and let $k$ and $l$ be real numbers such that
$k < l$. Then there exists a constant $\gamma >0$, depending only on $N$, such that
\begin{equation*}
(l-k) |K_1(y)\cap\{u\leqslant k\}|\,\, |K_1(y)\cap\{u\geqslant l\}| \leqslant \gamma  \int\limits_{K_1(y)\cap\{k\leqslant u\leqslant l\}}|\nabla u| dx.
\end{equation*}
\end{lemma}
\noindent For the next lemma we refer to \cite{CiaMosVes}, Lemma $4.1$.
\begin{lemma}\label{lem2.4}
For any $\beta>0$  there exists a constant $\gamma_0=\gamma_0(\beta)>1$  such that, for any bounded $u : Q_{1,1}(x_0,t_0)\rightarrow \mathbb{R}$ with $u(x_0, t_0)\geqslant 1$,  there exist $(y, \tau)\in Q_{1,1}(x_0, t_0)$  and $r > 0$ such that
$$Q_{r, r^{p}}(y, \tau)\subset Q_{1,1}(x_0, t_0),\qquad r^\beta \sup\limits_{Q_{r, r^{p}}(y, \tau)}u\leqslant \gamma_0\qquad
\mbox{and} \qquad r^\beta u(y, \tau)\geqslant \frac{1}{\gamma_0}.$$
\end{lemma}
\noindent In order to manipulate the nonlinearity in time, we will need also the following result which can be found \cite[Lemma 3.2]{BogDuzGiaLiaSch}

\begin{lemma}\label{lem2.5}
There exists a constant $\gamma >0$, depending only on $m$, such that
\begin{equation*}
\frac{1}{\gamma}\,(k^m+u^m)^{\frac{1}{m}-1}(u^m-k^m)^2_{\pm}\leqslant g_{\pm}(u^m, k^m)\leqslant \gamma\, (k^m+u^m)^{\frac{1}{m}-1}(u^m-k^m)^2_{\pm}
\end{equation*}
\end{lemma}
\noindent where
$$g_{\pm}(u^m, k^m):=\pm\,\frac{1}{m}\int\limits^{u^m}_{k^m} s^{\frac{1}{m}-1} (s-k^m)_{\pm}\,ds , \qquad k,m>0 \ .$$

\section*{Acknowledgements} \noindent S. Ciani acknowledges the partial funding of GNAMPA (INdAM) and the partial ones of the Math Department of the University of Bologna. E. Henriques was founded by Portuguese Funds through FCT - Funda\c c\~ao para a Ci\^encia e a Tecnologia - within the Projects UIDB/00013/2020 and UIDP/00013/2020 with DOI references 10.54499/UIDB/00013/2020 and 10.54499/UIDP/00013/2020. Igor Skrypnik and Mariia Savchenko were supported by the National Research Foundation of Ukraine, project No. $2025.07/0369$ "Qualitative methods of nonlinear analysis of heterogeneous structures". We thank professor M. Vestberg for precise suggestions regarding the existence theory underlying the operator under our study. Finally, we acknowledge the advice of the two anonymous referees, since their work has constructively helped us to improve ours.

\end{document}